\documentclass[winedt]{iitparc}
\usepackage{amssymb,amsmath}

\begin{document}

\frontmatter          
\IssuePrice{25}%
\TransYearOfIssue{2000}%
\TransCopyrightYear{2000}%
\OrigYearOfIssue{2000}%
\OrigCopyrightYear{2000}%
\TransVolumeNo{61}%
\TransIssueNo{9}%
\OrigIssueNo{9}%
\TransPartNo{1}

\mainmatter              

\Rubrika{STOCHASTIC SYSTEMS}
\CRubrika{STOCHASTIC SYSTEMS}

\setcounter{page}{1424}

\newcommand\propro[2]{\medskip\par\PPR{#1}{\rm #2}\par}
\newcommand\prothe[2]{\medskip\par\PTH{#1}{\rm #2}\par}
\newcommand{\PLE}[1]{{\bf Proof of Lemma~#1}}
\newcommand\prolem[2]{\medskip\par\PLE{#1}{\rm ~#2}\par}
\newcommand\theoremprimed[2]{\par\hskip5mm{\bfseries Theorem #1.\ \ }{\it #2}\par\bigskip}
\newcommand\coroll[2]{\par\hskip5mm{\bfseries Corollary #1.\ \ }{\itshape #2}\par\bigskip}
\newcommand\axiom[2] {\bigskip\par{\bfseries #1.\ } {\rm #2}\par\bigskip}
\newcommand\axiomt[2]{\par{\bfseries #1.\ }         {\rm #2}\par\bigskip}

\def\f{\varphi}                                       
\def\D{\Delta}                                        
\def\G{\Gamma}                                        
\def\si{\sigma}                                       
\def\h{H}                                             
\def\e{\varepsilon}                                   
\def\esm{\epsilon}                                    
\def\l{\ell}                                          
\def\FF{\mathop{\cal F}\nolimits}                     
\def\PP{\mathop{\cal P}\nolimits}                     
\def\GG{\mathop{\cal G}\nolimits}                     
\def\TT{\mathop{\cal T}\nolimits}                     
\def\LL{\mathop{\cal L}\nolimits}                     
\def\ul     {\mathop{\rm u}  \limits}                 
\def\ulo{\ul^{\scriptscriptstyle \circ}}              
\def\cupl   {\mathop{\cup}   \limits}                 
\def\maxl   {\mathop{\max}   \limits}                 
\def\suml   {\mathop{\sum}   \limits}                 
\def\capo { \mathop {\bigcap}\limits}                 
\def\cupo { \mathop {\bigcup}\limits}                 
\def\1n{1,\ldots,n}                                   
\def\_#1{\mathop{\hspace{-2pt}^{}_{#1}}}              
\def\S{\Sigma}                                        
\def\cdc{,\ldots,}                                    
\def\parr{\newline\indent}                            

\def\N{\mathop{\mathbb N}{}}                          
\def\R{\mathop{\mathbb R}{}}                          
\def\vj {\mathop{\bar{J}}\nolimits}                   
\def\q  {\mathop{\bar{J}}\nolimits}                   
\def\ktil {\widetilde K}                              
\def\intercal{\mathop{\scriptscriptstyle T}\nolimits} 

\raggedbottom

\author{R. P. Agaev \and P. Yu. Chebotarev}

\title{The Matrix of Maximum Out Forests of a Digraph and Its Applications%
\thanks{
This work was supported by the European Community under Grant INTAS--96--106.}}

\institute{Trapeznikov Institute of Control Sciences, Russian Academy of Sciences, Moscow, Russia}
\received{Received January 17, 2000}
\titlerunning{The Matrix of Maximum Out Forests of a Digraph}
\authorrunning{Agaev, Chebotarev}
\OrigCopyrightedAuthors{Agaev, Chebotarev} \OrigPages{15--43}


\maketitle

\begin{abstract}
{We study the maximum out forests of a (weighted) digraph and the matrix of maximum out forests. A maximum
out forest of a digraph $\G$ is a spanning subgraph of $\G$ that consists of disjoint diverging trees and has
the maximum possible number of arcs. If a digraph contains out arborescences, then maximum out forests
coincide with them. We consider {\it Markov chains related to a weighted digraph\/} and prove that the matrix
of Ces\`aro limiting probabilities of such a chain coincides with the normalized matrix of maximum  out
forests. This provides an interpretation for the matrix of Ces\`aro limiting probabilities of an arbitrary
stationary finite Markov chain in terms of the weight of maximum out forests. We discuss the applications of
the matrix of maximum out forests and its transposition, the {\it matrix of limiting accessibilities of a
digraph}, to the problems of preference aggregation, measuring the vertex proximity, and uncovering the
structure of a digraph.}
\end{abstract}

\section{Introduction}
\label{sec1}

The concept of maximum out forest of a digraph directly generalizes the notion of spanning diverging tree
(out arborescence), which is one of the central notions in the theory of directed graphs.  If spanning
diverging trees of a digraph exist, then they coincide with maximum out forests; otherwise maximum out
forests share their major properties. We study these properties in this paper, which has the following
structure. After the main notation, in Sections~\ref{forest} and~\ref{nata1} we study the properties of
spanning diverging forests, in Section~\ref{sec_constr} we give a block algorithm for their construction,
Section~\ref{nata} presents matrix-forest theorems, and in Section~\ref{qnc} we study the matrix of maximum
out forests. The main result of Section~\ref{sec7} states that the normalized matrix of maximum out forests
of a (weighted) digraph $\G$ coincides with the matrix of Ces\`aro limiting transition probabilities of any
Markov chain {\it related to\/}~$\G$. In Section~\ref{dosti}, the total weight of maximum out forests that
connect two vertices is considered as a measure of vertex accessibility. Sections~\ref{leader} and~\ref{stru}
deal with the applications of the matrix of maximum out forests in the contexts of scoring based on paired
comparisons and detecting the structure of digraphs.

\section{Notation}
\label{sec2}

\subsection{General terms}

\setcounter{footnote}{1} In the terminology, we mainly follow~\cite{Harary}. Suppose that $\G$ is a weighted
digraph without loops, $V(\G)=\{\1n\}$ $(n>1)$ is its set of vertices, and $E(\G)$ its set of arcs. The
weights of all arcs are supposed to be strictly positive. {\it A subgraph\/}\footnote{In the literature,
(see, e.g.,
\cite{Zykov69}) this object is sometimes called a {\em part\/} of $\G$, whereas a subgraph $\G'$ of digraph
$\G$ is defined as a part whose vertex set $V(\G')$ is a subset of $V(\G)$, whereas the arc set contains all
the arcs of $E(\G)$ that have both incident vertices belonging to $V(\G')$. Such a subgraph $\G'$ of digraph
$\G$ will be called here a {\it restriction\/} of $\G$ to $V(\G')$.} of a digraph $\G$ is a digraph whose
vertices and arcs respectively belong to the sets of vertices and arcs of~$\G$. A {\it spanning\/} subgraph
of $\G$ is a subgraph of $\G$ with vertex set $V(\G)$. The {\it indegree\/} id($w$) of vertex $w$ is the
number of arcs that come in~$w$. A~vertex $w$ will be called {\it undominated\/} if id($w$)=0 and {\it
dominated\/} if id$(w)\ge1$. A vertex $w$ is {\it isolated\/} in $\G$ if $\G$ does not contain arcs incident
to~$w$.

A {\it route\/} in a digraph is an alternating sequence of vertices and arcs $w_0, e_1,$ $w_1\cdc e_k, w_k$
with every arc $e_i$ being $(w_{i-1},w_i)$. A {\it path\/} in a digraph is a route all whose vertices are
different. A~{\it circuit\/} is a route with $w\_0=w_k$, the other vertices being distinct and different from
$w_0$. A~vertex $w$ {\it is reachable\/} from a vertex $z$ in $\G$ if $w=z$ or $\G$ contains a path from $z$
to~$w$. A {\it semipath} is an alternating sequence of distinct vertices and arcs, $w_0, e_1, w_1\cdc e_k,
w_k,$ where every arc $e_i$ is either $(w_{i-1},w_i)$ or $(w_i,w_{i-1}$). {\it Semicircuit\/} is defined in
the same way.

Let $E=(\e_{ij})$ be the matrix of arc weights. Its entry $\e_{ij}$ equals zero if and only if there is no
arc from vertex $i$ to vertex~$j$ in $\G$. If $\G^{\prime}$ is a subgraph of $\G$, then the weight of
$\G^{\prime}$, $\e(\G^{\prime})$, is the product of the weights of all its arcs; if $\G'$ contains vertices,
but does not contain arcs, then $\e(\G')=1$. The weight of a nonempty set of digraphs $\GG$ is defined as
follows:
\[
\e(\GG)=\suml_{H\in\GG}\e(H);
\]
the weight of the empty set is~0.

The {\it Kirchhoff matrix} \cite{Tutte} of a weighted digraph $\G$ is the $n\times n$-matrix
$L=L(\G)=(\l\_{ij})$ with elements $\l\_{ij}=-\e_{ji}$ when $j\ne i$ and $\l\_{ii}=- \suml_{k\ne i}\l\_{ik}$,
$i,j=\1n$.

\subsection{The structure of a digraph}

A {\em vertex basis\/} of a digraph $\G$ is any minimal (by inclusion) collection of vertices of $\G$ from
which all its vertices are reachable. The requirement of minimality can be equivalently replaced with that of
mutual unreachability of all vertices in the collection.

A digraph is called {\it strongly connected\/} (or {\it strong}) if all its vertices are mutually reachable,
{\em unilaterally connected\/} if for any two its vertices at least one of them is reachable from the other,
and {\it weakly connected\/} if any two different vertices are connected by a semipath.

The restriction of $\G$ to any equivalence class of the vertex mutual reachability relation is called a {\it
strong component}, or a {\it bicomponent}, or an {\it oriented leaf\/} of~$\G$. {\it Weak components\/} of
{$\G$} are defined similarly on the base of the vertex connectedness by semipaths. The unilateral
reachability relation can be intransitive, so it is not generally an equivalence relation. Nevertheless,
maximal (by the inclusion of vertex  sets) unilaterally connected subgraphs of $\G$ are sometimes called {\it
unilateral components\/} of $\G$. As distinct from strong and weak components, they may overlap.

Suppose that $\G_1\cdc\G_r$ are all the strong components of~$\G$. The {\it condensation\/} (or {\it
factorgraph}, or {\it leaf composition}, or {\it Hertz graph}) $\G^*$ of digraph $\G$ is the digraph with
vertex set $\{\G_1\cdc\G_r\}$ where an arc $(\G_i,\G_j)$ belongs to $E(\G^*)$ iff $E(\G)$ contains at least
one arc from a vertex of $\G_i$ to a vertex of~$\G_j$. The condensation of any digraph $\G$ contains no
circuits.

If a digraph does not contain circuits, then its vertex basis is obviously unique and coincides with the set
of all undominated vertices~\cite{Harary,Zykov69}. That is why the strong components of $\G$ that correspond
to undominated vertices of $\G^*$ are sometimes called the {\it basis bicomponents\/} of~$\G$~\cite{Zykov69}.
In this paper, the term {\it undominated knot of\/}~$\G$ will stand for the set of vertices of any basis
bicomponent of~$\G$:

\begin{definition}
\label{De1} A nonempty subset of vertices $K\subseteq V(\G)$ of digraph $\G$ is an {\em undominated knot\/}
in $\G$ if all the vertices that belong to $K$ are mutually reachable and there are no arcs $(w_j,w_i)$ with
$w_j\in V(\G)\setminus K$ and $w_i\in K$.
\end{definition}

An extreme case of undominated knot is a singleton consisting of an undominated vertex (if $\G$ contains such
vertices). The opposite extreme case is the whole vertex set of a strong digraph.

The following statement \cite{Harary,Zykov69} characterizes all the vertex bases of a digraph.

\begin{proposition}
\label{proZy} A set $W\subseteq V(\G)$ is a vertex basis of $\G$ if and only if\/ $W$\! contains exactly one
vertex from every undominated knot of $\G$ and no other vertices.
\end{proposition}

Schwartz \cite{Schwartz} refers to the undominated knots of a digraph as {\em minimum $P$-undominated sets}.
He formulates the Generalized Optimal Choice Axiom (GOCHA). If a preference relation (digraph) defined on a
finite set of alternatives is given, then the {\it choice\/} according to GOCHA is the union of minimum
$P$-undominated sets of this digraph.\footnote{This union is also called the {\em top cycle\/} and the {\em
strong basis\/} of the digraph.} This choice is interpreted as the set of ``best'' (in terms of GOCHA)
alternatives. A review of choice rules of this kind can be found in~\cite{VIV}.

\subsection{Diverging forests of a digraph}

A {\it diverging tree\/} is a digraph without semicircuits that has a vertex (called the {\it root}) from
which every its vertex is reachable. It is easy to see that the root is unique and its indegree is zero, the
indegrees of all other vertices being one. A~diverging tree is said to {\it diverge\/} from its root. A~{\it
diverging forest\/} is a digraph without circuits such that id$(w)\le1$ for every its vertex~$w$.

Let $F$ be a diverging forest. By indicating the vertices $w$ in $F$ such that id$(w)=0$ (these are called
the {\it roots of $F$}) and the subsets of vertices reachable from each root, we obtain a partition
$\{V_1(F)\cdc V_{v'}(F)\}$ of the vertex set $V(F)$ such that there exists a semipath in $F$ between $w\in
V_i(F)$ and $z\in V_j(F)$ if and only if $i=j$. Thus, the restriction of $F$ to every subset $V_i(F)$,
$i=1\cdc v'$, is a weak component of~$F$. It is easily seen that every component of a diverging forest is a
diverging tree.

For a fixed digraph $\G$, consider spanning diverging forests $F$ of~$\G$ (such subgraphs do obviously exist
for every digraph).

\begin{definition}
\label{De2} A spanning diverging forest $F$ of a digraph $\G$ is called a {\em maximum out forest\/} of $\G$
if $\G$ has no spanning diverging forest with a greater number of arcs than in~$F$.
\end{definition}

Obviously, every maximum out forest of $\G$ has the minimum possible number of roots; this number will be
called the {\it forest dimension\/}\footnote{This name recalls ``$W$-bases'', the term Fiedler and
Sedl\'{a}\v{c}ek \cite{Fiedler} used for spanning diverging forests.} of the digraph and denoted by~$v$. The
number of arcs in any maximum out forest is obviously $n-v$.

Let us emphasize that the property to be a maximum out fores is more stringent than the maximality with
respect to the inclusion of arc sets. This point is illustrated in the next section.

By $\FF(\G)=\FF$ and $\FF_k(\G)=\FF_k$ we will denote the sets of all spanning diverging forests of $\G$ and
the set of all spanning diverging forests of $\G$ with $k$ arcs, respectively; $\FF^{i\rightarrow j}_k$ will
designate the set of all spanning diverging forests with $k$ arcs where $j$ belongs to a tree diverging
from~$i$.

\section{Simple properties of diverging forests}
\label{forest}

\begin{lemma}
\label{lem1} Let $F$ be a diverging forest. $1.$ If a digraph $F'$ is obtained from $F$ by the removal of an
arc$,$ then $F'$ is a diverging forest too. $2.$ Suppose that digraph $F'$ is obtained from $F$ by the
addition of some arc $(z,w)$. In this case$,$ $F'$ is a diverging forest if and only if $w$ is undominated in
$F$ and $z$ is unreachable from~$w$.
\end{lemma}

The properties that make up Lemma~\ref{lem1} are obvious; we will use them in the proofs of Lemma~\ref{lem2}
and other statements without explicit references. The proofs are given in the Appendix.

\begin{lemma}
\label{lem2} If $w$ is a dominated vertex~of $\G,$ then for any $k\in\{\1n-v\},$ the set $\FF_k$ contains a
spanning diverging forest $F$ such that $w$ is dominated in $F$.
\end{lemma}

The stronger statement saying that for every $k$ and every arc $(z,w)\in E(\G)$, there exists a forest in
$\FF_k$ that contains $(z,w),$ is generally wrong. Indeed, consider the digraph shown in Fig.~1a. The forest
dimension of this digraph is one, and the unique maximum out forest $F$ (which is a diverging tree) is shown
in Fig.~1b. Arc $(4,2)$ is not in $F$. This demonstrates that for some spanning diverging forests $F'$, there
is no maximum out forest $F$ such that $E(F')\subseteq E(F)$. For instance, the arc sets of the spanning
diverging forests in Fig.~1c~--~1e are not contained in $E(F)$.

Thus, a maximal (with respect to the inclusion of arc sets) out forest can be not maximum (see, e.g.,
Fig.~1c,d). This implies, in particular, that the arc sets of spanning diverging forests of a digraph cannot
be considered as the independent sets of a matroid.

\addvspace{10mm}
\input rpe_f1e.pic

\hspace{70mm} Figure~1
\medskip

\begin{proposition}
\label{prop1} {\rm 1.} Any undominated vertex of a digraph is the root in every spanning diverging forest.
{\rm 2.} If a digraph does not contain circuits$,$ then no dominated vertex can be the root in a maximum out
forest.
\end{proposition}

\begin{lemma}
\label{lem3} {\rm 1.} The weights $\e(\FF_{k}^{t\rightarrow t})$ of the sets $\FF^{t\to t}_k$ are the same
for all undominated vertices $t$ of $\G$ and are equal to $\e(\FF_{k})$. {\rm 2.} For any $k\in\{\1n-v\},$
$\e(\FF_{k}^{t\rightarrow t})>\e(\FF_{k}^{w\rightarrow w})$ whenever $t$ is an undominated vertex in~$\G$ and
$w$ is a dominated vertex.
\end{lemma}

In the following statements, $i,j$, and $k$ are arbitrary vertices of~$\G$.

\begin{lemma}
\label{lem4} If there exists a path from $i$ to $j$ in $\G,$ and there is no path from $i$ to $j$ in a
maximum out forest $F$ of $\G,$ then for some $k \ne i,$ $F$ contains the arc $(k,j)$ or $i$ is reachable
from $j$ in~$F$.
\end{lemma}

Lemma~\ref{lem4} implies

\begin{proposition}
\label{prop2} If $i$ and $j$ belong to different trees in a maximum out forest $F$ of a digraph $\G$ and $j$
is a root in $F,$ then $\G$ contains no paths from $i$ to~$j$.
\end{proposition}

After replacing the hypothesis of Lemma~\ref{lem4} with its local version, we can give this lemma the form of
necessary and sufficient condition:

\begin{lemma}
\label{lem5} For any maximum out forest $F$ of digraph $\G$ and any vertices $i,j\in V(\G),$ $F$ does not
contain the arc $(i,j)$ that belongs to $E(\G)$ if and only if $F$ contains an arc $(k,j)$ for some $k\ne i$
or $i$ is reachable from $j$ in~$F$.
\end{lemma}

\section{Maximum out forests, bases, and undominated knots}
\label{nata1}

Suppose that $\ktil=\cupo^{u}_{i=1}K_i$, where $K_1\cdc K_u$ are all undominated knots of digraph~$\G$, and
$K_i^{+}$ is the set of all vertices reachable from $K_i$ and unreachable from the other undominated knots.
If $k\in\ktil$, then $K(k)$ will designate the undominated knot that contains~$k$. For any undominated knot
$K$ of $\G,$ denote by $\G_K$ the restriction of $\G$ to $K$ and by $\G_{-K}$ the subgraph with vertex set
$V(\G)$ and arc set $E(\G)\setminus E(\G_K)$. For a fixed $K$, $\TT$ will designate the set of all spanning
diverging trees of $\G_K$ and $\PP$ will be the set of all maximum out forests of $\G_{-K}$. By $\TT^k$
$(k\in K)$ we will denote the subset of $\TT$ consisting of all trees that diverge from $k$, and by $\PP^{K
\rightarrow i}$ ($i\in V(\G)$) the set of all maximum out forests of $\G_{-K}$ such that $i$ is reachable
from some vertex that belongs to~$K$ in these forests.

\begin{proposition}
\label{roots_basis} A set $W\subseteq V(\G)$ is the set of roots of a maximum out forest in $\G$ if and only
if $W$ is a vertex basis of~$\G$.
\end{proposition}

In view of Proposition~\ref{roots_basis}, the sets of roots of maximum out forests are characterized by
Proposition~\ref{proZy}. The following three statements follow from Propositions~\ref{proZy}
and~\ref{roots_basis}.

\begin{proposition}
\label{prop6a} For any maximum out forest $F$ of a digraph $\G$ and any undominated knot $K_i,$ the
restriction of $F$ to $K^{+}_i$ is a diverging tree.
\end{proposition}

\begin{proposition}
\label{prop4a} The forest dimension of a digraph is equal to the number of its undominated knots$:$ $v=u$.
\end{proposition}

\begin{proposition}
\label{prop5a}
The forest dimension of a strong digraph is one.
\end{proposition}

To prove Proposition~\ref{prop5a}, it suffices to observe that the unique undominated knot of a strong
digraph $\G$ is $V(\G)$.

Since every weak component can be split into strong components, at least one of which being an undominated
knot, the forest dimension of a digraph is nonstrictly between the number of weak components and number of
strong components. The number of unilateral components is also no less than the number of weak components
(because every weak component contains at least one unilateral component), but it can be less, or greater, or
equal to the number of strong components; this number can even exceed the number of vertices in $\G$ (an
example is the bipartite digraph $\G$ where $V(\G)=\{i_1,i_2,i_3,j_1,j_2,j_3\}$ and $E(\G)=\{(i_k,j_t)\mid
k,t=1\cdc3\}$). Finally, the vertex set of any unilateral component either has the empty meet with $\ktil$ or
contains exactly one undominated knot. That is why the forest dimension of a digraph cannot exceed the number
of its unilateral components. Thus, the following statement holds true.

\begin{proposition}
\label{propner} The forest dimension of a digraph is no less than its number of weak components and does not
exceed the number of its strong components and the number of its unilateral components.
\end{proposition}

Let us fix an arbitrary undominated knot $K$ of $\G$ and consider the sets $\TT,\PP,\TT^k,$ and $\PP^{k\to
i}$ defined above.

Let $\TT\odot\PP=\{T\cup F: T\in\TT, F\in\PP\}$, where $T\cup F$ is a digraph with vertex set $V(T)\cup
V(F)=V(\G)$ and arc set $E(T)\cup E(F)$. In the same way, $\TT^k\odot\PP^{K\rightarrow i}=\{T\cup F:
T\in\TT^k, F\in\PP^{K\rightarrow i}\}$.

\begin{proposition}
\label{prop4} Suppose that $K$ is an arbitrary undominated knot of $\G$ and the sets $\TT,$ $\PP,$ $\TT^j,$
and $\PP^{K\rightarrow i}$ $(j\in K,$ $i\in V(\G))$ are determined by~$K$. Then

{\rm 1.} $\FF_{n-v}=\TT\odot\PP$ and $\e(\FF_{n-v})=\e(\TT)\e(\PP);$

{\rm 2.} for any $j\in K$ and $i\in V(\G),$ we have
\begin{equation}
\label{compos} \FF^{j\rightarrow i}_{n-v}=\TT^{j}\odot\PP^{K \rightarrow i}\quad\mbox{\rm and}\quad
\e(\FF_{n-v}^{j\rightarrow i})=\e(\TT^j)\e(\PP^{k\rightarrow i}).
\end{equation}
\end{proposition}

\section{An algorithmic description of maximum out forests}
\label{sec_constr}

In this section, we give a block algorithm for constructing all maximum out forests of~$\G$.

1.~Find all undominated knots $K_1\cdc K_v$ of $\G$ and the sets $K_1^+\cdc K_v^+$.

2.~In every $K_i^+,$ construct an arbitrary spanning diverging tree (rooted within $K_i$).

3.~Find the strong components in the restriction of $\G$ to $V(\G)\setminus\,\cupl_{i=1}^v\,K_i^+$. Let
$T_1\cdc T_s$ be the sets of vertices of these strong components.

4.~For every $T_i,$ draw one or any greater number of arcs taken from $E(\G)$ and directed to {\it
distinct\/} vertices in~$T_i$.

5.~For every $T_i,$ construct an arbitrary spanning forest rooted at those and only those vertices to which
the arcs from outside were drawn on step~4.

6.~Consider the spanning subgraph $F$ whose arc set consists of all arcs drawn on steps~2, 4 and~5.

Step~5 can be reduced to the construction of a tree in the following way.

5a.~Identify all the vertices of $T_i$ to which the arcs from outside were drawn on step~4. Let the resulting
vertex be~$t^*_i$. Construct an arbitrary diverging tree spanning in the remaining part of $T_i$ and rooted
at~$t^*_i$. Now split $t^*_i$ into the vertices constituting it and replace the arcs directed from $t^*_i$
with arbitrary corresponding arcs directed from these vertices.

\begin{proposition}
\label{algo} $1.$~The sets of subgraphs determined by steps~$5$ and $5${\rm{a}} of the above block algorithm
coincide.

$2.$~The set of subgraphs produced by the block algorithm $1$--$6$ coincides with set of maximum out forests
of~$\G$.
\end{proposition}

\section{Parametric versions of the matrix-forest theorem}
\label{nata}

In \cite{CheSha981} we presented a parametric version of the matrix-forest theorem for multigraphs:
\begin{theorem}
\label{th1} For any weighted multigraph~$G$ with positive weights of edges and any $\tau>0,$ there exists the
matrix $Q(\tau)=(q_{ij}(\tau))=(I+{\tau}L(G))^{-1}$ and
\[
q_{ij}(\tau)=\sum^{n-v}_{k=0}\tau^k{\e}(\FF^{ij}_k) \Big/\sum^{n-v}_{k=0}\tau^k{\e}(\FF_k),\quad i,j=\1n,
\]
where $\FF_k$ is the set of all spanning rooted forests of $G$ that contain $k$ edges$,$ $\FF^{ij}_k$ is the
set of all $k$-edge spanning rooted forests of $G$ where $j$ belongs to a tree rooted at $i,$ and $v$ is the
number of components in~$G$.
\end{theorem}

An analogous theorem is true for {\it multidigraphs\/} (that may contain multiple arcs between different
vertices, but not loops).
\bigskip

{\bf Theorem~\ref{th1}$'$.} {\it For any weighted multidigraph $\G$ with positive weights of arcs and any
$\tau>0,$ there exists the matrix $Q(\tau)=(q_{ij}(\tau))=(I+{\tau}L(\G))^{-1}$ and
\begin{equation}
q_{ij}(\tau)=\sum^{n-v}_{k=0}\tau^k{\e}(\FF^{j\rightarrow i}_k) \Big/\sum^{n-v}_{k=0}\tau^k{\e}(\FF_k),\quad
i,j=\1n, \label{param}
\end{equation}
where $\FF_k$ and $\FF^{j\rightarrow i}_k$ are defined at the end of Section~$\ref{sec2},$ and $v$ is the
forest dimension of\/~$\G$.}
\bigskip

To prove this theorem, it suffices to apply the matrix-forest theorem for multidigraphs \cite{CheSha97} to
the weighted multidigraph $\G'$ that differs from $\G$ in the weights of arcs only: for all $i,j=\1n,$
$\e'_{ij}=\tau\e\_{ij}$.

The matrix form of this theorem is as follows:
\bigskip

{\bf Theorem~\ref{th1}$''$.} {\it For any weighted multidigraph $\G$ with positive weights of arcs and any
$\tau>0,$ there exists the matrix $Q(\tau)=(I+{\tau}L(\G))^{-1}$ and
\[
Q(\tau)=\frac{1}{s(\tau)}\left({\tau}^0 Q_0 +{\tau}^1 Q_1+\ldots+{\tau}^{n-v} Q_{n-v}\right), \label{razlo}
\]
where
\begin{equation}
\label{raz2} s(\tau)=\suml^{n-v}_{k=0}{\tau}^k\e(\FF_k), \;  Q_k=(q^k_{ij}), \;
q^k_{ij}=\e(\FF^{j\rightarrow i} _k),\; k=0\cdc n-v,\>
    i,j=\1n,
\end{equation}
and $\FF_k$ and $\FF_k^{j\rightarrow i}$ are the same as in Theorem~{\rm\ref{th1}}$'$. }
\bigskip

In the case of undirected graphs, the entries of the matrix of maximum rooted forests $Q_{n-v}$ are the same
within every component of~$G$. In the directed case, the matrix $Q_{n-v}$ possesses nontrivial properties
determined by the properties of maximum out forests. This matrix is studied in the following three sections.

\section{The matrix of maximum out forests}
\label{qnc}

According to (\ref{raz2}), $Q_{n-v}=(q^{n-v}_{ij})$, where $q^{n-v}_{ij}=\e(\FF^{j\rightarrow i}_{n-v})$,
i.e., the element $q_{ij}^{n-v}$ of $Q_{n-v}$ is the weight of the set of all {\it maximum\/} out forests of
digraph~$\G$ such that $i$ belongs to a tree diverging from~$j$. That is why $Q_{n-v}$ can be called the {\it
matrix of maximum out forests\/} of~$\G$.

\begin{theorem}
\label{th2} Suppose that $\G$ is an arbitrary digraph and $K$ is an undominated knot in~$\G$. Then the
following statements are true:

{\rm 1.}~For any $i\in V(\G),\;$ $\suml^n_{j=1}q^{n-v}_{ij}=\e(\FF_{n-v}).$

{\rm 2.}~$q^{n-v}_{ij}\ne 0\;\Leftrightarrow\; (j\in \ktil $ and $i$ is reachable from $j$ in $\G).$

{\rm 3.}~Suppose that $j\in K.$ Then for any $i\in V(\G),$
 $q^{n-v}_{ij}=\e(\TT^j)\e(\PP^{K\rightarrow i})$. Moreover$,$ if
$i\in K^{+},$ then $q^{n-v}_{ij}=q^{n-v}_{jj}=\e(\TT^j)\e(\PP).$

{\rm 4.}~$\suml_{j\in K}q^{n-v}_{jj}=\e(\FF_{n-v})$. In particular$,$ if $j$ is undominated$,$ then
$q^{n-v}_{jj}=\e(\FF_{n-v}).$

{\rm 5.}~If $j_1,j_2\in K,$ then $q^{n-v}_{\cdot j\_2}=(\e(\TT^{j_2}) \slash\e(\TT^{j_1}))q^{n-v}_{\cdot
j\_1},$ i.e.$,$ the $j_1$ and $j_2$ columns of $Q_{n-v}$ are proportional{\rm.}
\end{theorem}

Note that if the forest dimension of $\G$ is~1, i.e., $\G$ contains a spanning diverging tree, then
$Q_{n-v}=Q_{n-1}$ and $q^{n-1}_{ki}=q^{n-1}_{ji}$ for all $i,j,k\in V(\G)$. Indeed, in this case,
$q^{n-1}_{ji}$ is a the total weight $\e(\TT^i)$ of all spanning trees diverging from~$i$. Therefore, by the
matrix-tree theorem, $Q_{n-1}$ coincides in this case with the matrix of cofactors (the adjugate matrix)
of~$L$.

\begin{definition}
\label{De3} The matrix $\vj=(\q_{ij})=\si^{-1}Q_{n-v},$ where $\si=\e(\FF_{n-v}),$ will be called the {\it
normalized matrix of maximum out forests\/} of a digraph.
\end{definition}

The matrix $\vj$ will be the focus of our attention in what follows. First of all, we reformulate
Theorem~\ref{th2} for~$\vj$.
\bigskip

\theoremprimed{\ref{th2}$'$} {Suppose that $\G$ is an arbitrary digraph and $K$ is an undominated knot
in~$\G$. Then the following statements are true.

{\rm 1.}~$\vj$ is a stochastic matrix$:$ $\q_{ij}\ge0,$ $\;\suml^n_{k=1}\q_{ik}=1,\;$ $i,j=\1n.$

{\rm 2.}~$\q_{ij}\ne 0\;\Leftrightarrow\; (j\in \ktil$ and $i$ is reachable from $j$ in~$\G).$

{\rm 3.}~Suppose that $j\in K.$ For any $i\in V(\G),$
 $\q_{ij}=\e(\TT^j)\e(\PP^{K\rightarrow i})\slash
\e(\FF_{n-v})$. Furthermore$,$ if $i\in K^{+},$ then $\q_{ij}=\q_{jj}=\e(\TT^j)\slash \e(\TT)${\rm.}

{\rm 4.}~$\suml_{j\in K}\q_{jj}=1.$ In particular$,$ if $j$ is an undominated vertex$,$ then $\q_{jj}=1.$

{\rm 5.}~If $j_1,j_2\in K$, then $\q_{\cdot j\_2}=(\e(\TT^{j_2}) \slash\e(\TT^{j_1}))\q_{\cdot j_1},$ i.e.$,$
the $j_1$ and $j_2$ columns of $\q$ are proportional{\rm.}}
\medskip

Theorem~\ref{th2}$'$ follows from Theorem~\ref{th2}. To prove the last statements of item~3, item~1 of
Proposition~\ref{prop4} can be additionally used.
\bigskip

\begin{coroll}
{\bf from item~3 of Theorem~\ref{th2}$'$ and Proposition~\ref{prop6a}} {{\rm 1.} The normalized matrix of
maximum out forests $\vj^K=(\vj^K_{ij})$ of $\G_K$ coincides with the principal submatrix of $\vj$
corresponding to~$K.$

{\rm 2.} If $i\in K^+$ and $j\in K^+\setminus\,K,$ then $\vj$ is preserved under any variation of the weight
of $(i,j)$. }
\end{coroll}

Let $K(i)$ be the undominated knot that includes~$i$, provided that $i\in\ktil$. The following theorem is
concerned with the comparison of the entries of~$\vj$.

\begin{theorem}
\label{th3} For any $\G$ and any $i,j\in\{\1n\},$ the following statements are true.

{\rm 1.} $\q_{ii}\ge\q_{ji}.$

{\rm 2.} If $\q_{ii}>\q_{ji},$ then $i\in\ktil$ and $j\notin K^+(i),$ therefore$,$ $\G$ contains no paths
from $j$ to~$i.$

{\rm 3.} If $\q_{ii}>\q_{ji}>0,$ then $j\notin\ktil,$ consequently$,$ $j$ is not the root in any maximum out
forest of~$\G.$

{\rm 4.} If $\q_{ij}>0,$ then $\q_{ii}=\q_{ji}.$
\end{theorem}

\begin{theorem}
\label{I} For every weighted digraph$,$ $\vj$ is idempotent$:$ $\;\q^{2}=\q.$
\end{theorem}

Recall that $L=L(\G)=(\l\_{ij})$ is the Kirchhoff matrix of~$\G$.

\begin{theorem}
\label{th4} For every weighted digraph$,$ $L\q=\q L=0$.
\end{theorem}

It is worth noting a certain duality between $L$ and~$\vj$.

\begin{proposition}
\label{prop11} The ranks of $L$ and $\vj$ are $n-v$ and $v,$ respectively.
\end{proposition}

\begin{rmrk}
Consider the rows  ${\l\_1\cdc\l\_n}$ of $L$ as vectors in~$\R^n$. Denote by $\LL$ the multiset of these rows
and by $L_R$ the linear span of ${\l\_1\cdc\l\_n}$ in~$\R^n$. Since $\LL$ contains $n-v$ linearly independent
vectors (Proposition~\ref{prop11}), the dimension of $L_R$ is~$n-v$.

Let ${\vj}_R$ be the linear span of the columns of~$\vj.$ By Proposition~\ref{prop11}, the dimension of
${\vj}_R$ is~$v$.

Note that: (A) $L_R\cap {\vj}_R=\{0\}$. Indeed, if the meet of these two subspaces contained a nonzero vector
$u$, then, by Theorem~\ref{th4}, $\| u\|^2=0$ would hold; (B) the dimensions of $L_R$ and ${\vj}_R$ sum
to~$n$.

By (A) and (B), $\R^n$ is decomposable into the {\it direct sum\/} of the subspaces $L_R$ and ${\vj}_R$ (see,
e.g.,~\cite{Gelf}):
\[
\R^n=L_R\dot+{\vj}_R,
\]
i.e., every vector $u\in\R^n$ can be uniquely represented as $u=u_1+u_2$, where $u_1\in L_R$ and $u_2\in\vj$.
\end{rmrk}

The following theorem provides an explicit expression for~$\vj$.

\begin{theorem}
\label{Tl} For any weighted multidigraph $\G,$
\begin{equation}
\label{l} \vj=\lim_{\tau\to\infty} (I+\tau\,L)^{-1}.
\end{equation}
\end{theorem}

The stochasticity and idempotence proven for $\vj$ are typical of the limiting transition probability
matrices of Markov chains. These matrices also possess some properties that resemble Theorem~\ref{th4} and
the other above statements. This resemblance is not accidental. It turns out that $\vj$ determines the
asymptotic behavior of certain Markov chains related to~$\G$. The corresponding results are presented in the
following section.

\section{Markov chains related to a weighted digraph}
\label{sec7}

\begin{definition}
\label{3} Let us say that a stationary Markov chain with set of states $\{\1n\}$ and transition probability
matrix $P$ is {\it related to a weighted digraph} $\G$ if there exists $\alpha\ne0$ such that
\begin{equation}
\label{7.1} P=I-\alpha\,L(\G).
\end{equation}
\end{definition}

We will identify the states of this Markov chain with the corresponding vertices of\/~$\G$. According to
Definition~\ref{3}, if a Markov chain is related to a weighted digraph $\G$, then the probability of
transition from $i$ to $j$ is proportional to the weight of the $(j,i)$ arc in~$\G$. Thus, if the weight of
the arc $(j,i)$ is interpreted as the degree of preference given to vertex $j$ in a comparison with $i$ or
something like that, then this weight determines the probability of transition from the {\it dominated\/}
vertex to the {\it dominating\/} one (the transitions are laid from the ``worse'' to the ``better'').

It is easy to see that the union of undominated knots of $\G$ is the set of {\it essential states\/} (in
Kolmogorov's notation) of any Markov chain related to~$\G$. All other vertices are {\it unessential states\/}
of every such a chain.

According to~(\ref{7.1}), the {\it row defect\/} $1-\alpha\sum_{j=1}^n\,\e_{ji}$ determines the probability
of transition from $i$ to $i$ (a stagnant transition).

Since we consider finite and stationary Markov chain only, we will omit the words ``finite'' and
``stationary.''

Definition~\ref{3} differs from the customary way of attaching Markov chains to graphs (used in
\cite[Chapter~9]{Bollo}, \cite{SearyRich}, and many other works). In the Markov chain attached to a graph in
accordance with the classical definition, the transition probabilities for the pairs of different vertices
are not generally proportional to the corresponding arc (edge) weights.
As a result, the transition probability matrices of the Markov chains attached to the graphs with symmetric
matrices $E$ and $L$ are generally nonsymmetric.

It is easy to determine the condition under which the matrix~(\ref{7.1}) represents the transition
probabilities of some Markov chain.

\begin{proposition}
\label{p7.1} Matrix $P$ defined by $(\ref{7.1})$ is the transition probability matrix of a Markov chain
$($and thus this chain is related to $\G$ in terms of Definition~$\ref{3})$ if and only if
$0<\alpha<(\maxl_{1\le\,i\le\,n}\,\l\_{ii})^{-1}$.
\end{proposition}

Proposition~\ref{p7.1} provides a necessary and sufficient condition for the stochasticity of
$P=I-\alpha\,L(\G)$; it immediately follows from the definition of~$L$.

\begin{rmrk}
Consider the matrix norm $\|\cdot\|_{\infty}$ on the set of $n\times n$ Kirchhoff matrices $L$:
\[
\|\,L\,\|_{\infty}=\maxl_{1\le i\le n}\,\suml_{j=1}^n\l\_{ij}= 2\maxl_{1\le i\le n}\l\_{ii}.
\]
This matrix norm is called the {\it maximum row sum norm\/}~\cite{HoJo}.

The function $\|\,L\,\|_{\omega}=\maxl_{1\le i\le n}\l\_{ii}$ appearing in Proposition~\ref{p7.1} is not a
matrix norm, since it does not obey submultiplicativity $\|\,AB\,\|\le\|\,A\,\|\|\,B\,\|$. Indeed, to
demonstrate this, it is sufficient to take for $A$ and $B$ the Kirchhoff matrix of the digraph on two
vertices with two symmetric arcs carrying unit weights. The other axioms of matrix norms (nonnegativity,
positivity, homogeneity, and triangle inequality) are satisfied for $\|\cdot\|_{\omega}$. Thus, this function
is a {\it generalized matrix norm\/}~\cite{HoJo}.
\end{rmrk}

\begin{proposition}
\label{Q1} The matrix $(\suml_{k,t=1}^n\e_{kt})^{-1} Q_1,$ where $Q_1$ is the matrix of diverging forests
with one arc $($defined in Theorem~$\ref{th1}''),$ is the transition probability matrix of some Markov chain
related to~$\G$.
\end{proposition}

It is easily seen that every Markov chain is related to some weighted digraph. More exactly, there is always
a family of such digraphs~$\G$: their Kirchhoff matrices are
\begin{equation}
\label{7.2} L(\G)=\frac{1}{\alpha}(I-P)
\end{equation}
\noindent with all possible $\alpha>0$. The matrices of arc weights of all these digraphs are proportional.

The sequence $P,\,P^2,\,P^3,\ldots$ is nonconvergent for periodic Markov chains. Consider the Ces\`aro limit
of this sequence, which can be shown to exist for every Markov chain.

\begin{definition}
\label{4} The {\it limiting matrix of average probabilities\/} of a Markov chain is the matrix
\begin{equation}
\label{B} B=\lim_{k\to\infty}\frac{1}{k}\,\suml_{t=0}^{k-1} P^t,
\end{equation}
provided that this limit exists.
\end{definition}

Consider also the matrix
\begin{equation}
\label{B1} B'=\frac{1}{m}\,\suml_{j=0}^{m-1} B_j,
\end{equation}
\noindent where $m$ is the period of the Markov chain and $B_0\cdc B_{m-1}$ are the limiting matrices for the
convergent subsequences of $\{P^t\}$:
\begin{equation}
\label{Bj} B_j=\lim_{i\to\infty} P^{im+j},\quad j=0\cdc m-1.
\end{equation}
The case $m=1$ corresponds to convergent sequences $\{P^t\}$.

The following statement is known.
\begin{proposition}
\label{p7.2} For every Markov chain$,$ there exists the limiting matrix of average probabilities $B,$ and
$B=B'$.
\end{proposition}

In the Appendix, we give a proof of this proposition which is closely connected with the proof of the main
result of this section (Theorem~\ref{M} below).
\bigskip

\begin{coroll}
{\bf 1 from Proposition~\ref{p7.2}} {If the sequence $\{P^t\}$ converges and $P^{\infty}$ is its limit$,$
then $P^{\infty}=B$.}
\end{coroll}

\begin{coroll}
{\bf 2 from Proposition~\ref{p7.2}} {For any Markov chain$,$

{\rm (1)} $PB=BP=B${\rm;}

{\rm (2)} the nonzero columns of $B$ are right eigenvectors$,$ the rows being left eigenvectors of $P,$ all
corresponding to eigenvalue~$1${\rm;}

{\rm (3)} for every weighted digraph $\G$ to which a Markov chain is related$,$ $BL(\G)=L(\G)B=0$ holds{\rm;}

{\rm (4)} $B$ is idempotent$:$ $B^2=B.$}
\end{coroll}

Let
\begin{equation}
\label{B(k)} B(k)=\bigl(b_{ij}(k)\bigr)=\frac{1}{k}\,\suml_{t=0}^{k-1} P^t,\quad k=1,2,\ldots
\end{equation}

A sort of experiment can be indicated where $B(k)$ and $B$ are the transition probability matrices. The
notion ``point in time'' means in the following statement the number of transition occurred in a Markov
chain.

\begin{proposition}
\label{7.p} {\rm 1.} Every element $b_{ij}(k)$ of $B(k)$ is the probability that the state of the Markov
chain at a random point in time uniformly distributed on $\{0,1\cdc k-1\}$ is $j,$ provided that the initial
state is~$i.$

{\rm 2.} The probability specified in item~{\rm 1} of this proposition tends to the element $b_{ij}$ of $B$
as $k\to\infty$.
\end{proposition}

Item~2 of Proposition~\ref{7.p} refers to experiments where the maximum possible number of Markov chain's
transitions antecedent to the instant of observation is not bounded {\em a priori}; this is a surrogate of
the impossible uniform distribution on a denumerable set.

\begin{proposition}
\label{7.3} For every Markov chain$,$ we have
\begin{equation}
\label{Bl} B=\lim_{\tau\to\infty}\bigl(I-\tau(P-I)\bigr)^{-1}.
\end{equation}
\end{proposition}

We are now in position to formulate the main result of this section.

\begin{theorem}
\label{M} \ For any Markov chain related to a weighted digraph $\G,$ the limiting matrix of average
probabilities $B$ coincides with~$\vj.$
\end{theorem}

Theorem~\ref{M} provides a method for a finite (combinatorial) calculation of $B$ (and thus of the stationary
distributions of Markov chains and of $P^\infty$, provided that the latter matrix exists). This method
consists in finding and classifying the maximum out forests of a digraph $\G$, i.e., in calculating~$\vj$.

The following corollary presents the fact that this method is applicable to every finite Markov chain.
\bigskip

\begin{coroll}
{\bf from Theorem~\ref{M}} {For any Markov chain$,$ the limiting matrix of average probabilities $B$ is equal
to the matrix $\vj$ of any weighted digraph to which this chain is related$,$ i.e.$,$ of any weighted digraph
$\G$ that has the Kirchhoff matrix $L(\G)=\alpha^{-1}(I-P),$ where $\alpha>0$ and  $P$ is the transition
probability matrix of the  Markov chain.}
\end{coroll}

This corollary is deduced from Theorem~\ref{M} and Definition~\ref{3}.
\smallskip

\begin{rmrk}
By virtue of the above corollary, Theorem~\ref{M}, and Proposition~\ref{p7.2}, $\vj$ can be substituted for
$B$ and $B'$ in all statements of this section. In this way, items~3 and~4 of Corollary~2 provide
Theorem~\ref{th4} and Theorem~\ref{I}, respectively. Thus, the Markov chain technique enables one to get
alternative proofs of these theorems.
\end{rmrk}

\section{The weight of maximum out forests as a measure of
vertex accessibility} \label{dosti}

By Theorem~\ref{M}, the matrix $\q=(\q_{ij})$ of a weighted digraph $\G$ coincides with the limiting matrix
of average probabilities of any Markov chain related to~$\G$. That is why $\q_{ij}$ can be called the {\it
limiting accessibility\/} of $i$ from $j$ in random walks on $\G$ with transition probabilities proportional
to the arc weights. The matrix $\q^{\intercal}$ will be referred to as the {\it matrix of limiting
accessibilities\/} in $\G$. In this section, we consider the entries of $\q^{\intercal}$ as a measure of
``proximity'' between vertices. For this purpose, we turn to the conditions proposed in
\cite{CheSha981,CheSha97} for the description of the notion of vertex proximity. These conditions are not
considered as necessary attributes of proximity measures, but if some index breaks a majority of them, this
indicates that the index measures not proximity, but something different.

\axiom{Nonnegativity} {For any digraph $\G,$ $p\_{ij}\ge0,\;\:i,j\in\{\1n\}$.}

\axiomt{Reversal property} {For any digraph $\G,$ the reversal of all its arcs $($provided that their weights
are preserved$)$ results in the transposition of the proximity matrix.}

\axiomt{Diagonal maximality} {For any digraph $\G$ and any distinct $i,j\in V(\G),$ $p\_{ii}>p\_{ij}$ holds.}

\axiomt{Triangle inequality for proximities} {For any digraph $\G$ and any $i,j,k\in V(\G)$,
$p\_{ij}+p\_{ik}-p\_{jk}$ $\le p\_{ii}$ holds. If, in addition, $j=k$ and $i\ne j$, then the inequality is
strict.}

Let
\[
d\_{ij}=p\_{ii}+p\_{jj}-p\_{ij}-p\_{ji},\quad i,j\in\{\1n\}.
\]

\axiomt{Metric representability of proximity} {The index $d\_{ij}$ is a distance between the vertices of a
digraph, i.e., it satisfies the axioms of metrics.}

It has been shown in \cite{CheSha982} that the triangle inequality for proximities corresponds to the
ordinary triangle inequality for the values~$d\_{ij}$.

\axiom{Disconnection condition} {For any digraph $\G$ and any $i,j\in V(\G),\;$ $p\_{ij}=0$ if and only if
$j$ is unreachable from~$i$.}

\axiomt{Transit property} {For any digraph $\G$ and any $i,k,t\in V(G),$ if $\G$ contains a path from $i$ to
$k,$ $i\ne k\ne t,$ and every path from $i$ to $t$ includes $k,$ then $p\_{ik}>p\_{it}.$}

The following condition is formulated here in a weaker version compared to that in~\cite{CheSha981}.

\axiom{Monotonicity} {Suppose that the weight of some arc $\e_{kt}^p$ in a digraph $\G$ increases. Then{\rm:}
\parr
{\rm 1)} $\D p\_{kt}>0,$ and for any $i,j\in\{\1n\},$ $(i,j)\ne (k,t)$ implies $\D p\_{kt}>\D p\_{ij};$
\parr
{\rm 2)} for any $i\in\{\1n\},$ if there is a path from $i$ to $k,$ and each path from $i$ to $t$ includes
$k,$ then $\D p\_{it}>\D p\_{ik}.$ }

Suppose that $P=(p_{ij})=\vj^{\intercal}$ is the matrix of limiting accessibilities of a digraph.

\begin{proposition}
\label{bliz} The index of limiting accessibilities of a digraph satisfies nonnegativity and the\/ {\rm
`$\Leftarrow$'} part of disconnection condition$;$ diagonal maximality$,$ transit condition$,$ and the first
part of item~$1$ and item~$2$ of monotonicity are satisfied in the nonstrict form$;$ reversal property$,$
triangle inequality for proximities$,$ metric representability of proximity$,$ the\/ {\rm `$\Rightarrow$'}
part of disconnection condition$,$ and the second part of item~$1$ of monotonicity are not satisfied.
\end{proposition}

In view of Proposition~\ref{bliz}, the index of limiting accessibilities does not completely correspond to
the concept of proximity lying in the above conditions. This is because it expresses accessibility {\em in
infinite time}. As we are going to show elsewhere, the replacement of
$\vj=\lim_{\tau\to\infty}(I+\tau\,L)^{-1}$ (Theorem~\ref{Tl}) by $(I+\tau\,L)^{-1}$ with a finite positive
$\tau$ results in a more sensible index of vertex proximity.

\section{The matrix of limiting accessibilities and the problem of determining leaders}
\label{leader}

Ranking players on the base of tournaments or irregular pairwise contests is an old, but still intriguing
problem. A statistical version of this problems is estimating objects on the base of paired
comparisons~\cite{David}. Analogous problems of the analysis of individual and collective preferences arise
in the contexts of voting, expert judgment, sociology, and psychometrics. Hundreds of methods have been
proposed for the solution of these problems (see, e.g.,
\cite{David,DavidAndr,CookKress,BelkLev,CheShaSpr,CheSha99,Laslier}).

In this paper, we suppose that an incomplete tournament with weighted results of paired contests or an
incomplete structure of numerical preferences is represented by a weighted digraph~$\G$.

One of the most popular sensitive methods for assigning scores to the participants in a tournament was
proposed by Daniels in 1969 and reduces, in our notation, to finding nonzero and nonnegative solutions to the
system of equations
\begin{equation}
L^{\scriptscriptstyle T} x=0. \label{D}
\end{equation}

The entry $x_i$ of the solution vector $x$ is used as an evaluation attached to the object represented by
vertex~$i$.

The system of equations (\ref{D}) in a componentwise notation has the form:
\begin{equation}
\l\_{ii} x_i=\suml_{j\ne i}\,(-\l\_{ji})x_j,\quad i=\1n. \label{eM}
\end{equation}

In the interpretation by Moon and Pullman \cite{MoonPullman}, $x_i$ is the amount player $i$ pays to any
player that defeats~$i$. In the simplest case of a nonweighted (but generally incomplete) tournament, the
left-hand side of (\ref{eM}) is the amount paid by the player for her defeats, whereas the right-hand side is
the amount collected by the player for her wins. Thus, the equality of these amounts for all participants
stated by (\ref{eM}) can be considered as a fairness condition\footnote{More exactly, this condition means
that the payoff vector is representative of the strength of the players.} imposed on the \textit{payoff
vector} $x=(x\_1\cdc x\_n)$: if the strength of each player remains the same, then nobody receives any
advantage, and everyone can expect a zero total.

This method was rediscovered several times with different motivations (some references are given
in~\cite{CheSha99}). As was noticed by Berman \cite{Berman} (although, in other contexts, this had been
remarked by Maxwell \cite{Maxwell} and other writers), if a tournament is strong, i.e. all its vertices are
mutually reachable, then the general solution to (\ref{D}) is given by the vectors proportional to
$t=(t\_1\cdc t\_n),$ where $t\_{i}$ is the weight of the set of spanning trees (out arborescences) diverging
from~$i$. This fact can be easily proved as follows. By the matrix-tree theorem for digraphs (see,
e.g.,~\cite{Harary}), $t\_j$ is the cofactor of any entry in the $j$th row of~$L$. Then for every $i\in
V(\G),$ $\suml^n_{j=1}\,\l\_{ij}\,t_j={\rm det}\,L$ (the expansion of the determinant by the $i$th column
of~$L$) and, since ${\rm det}\,L=0,\;$ $t$ is a solution to~(\ref{D}). As ${\rm rank}\,L=n-1$ (since the
cofactors of $L$ are nonzero), any solution to (\ref{D}) is proportional to~$t$.

Berman \cite{Berman} asserted that this result is sufficient to rank order the vertices in any digraph,
because its strong components supposedly ``can be ranked such that every player in a component of higher rank
defeats every player in a component of lower rank. Now by ranking the players in each component we obtain a
ranking of all the players.''

While the statement of the existence of a natural ranking of the strong components is correct in the case of
round-robin tournaments, it is obviously mistaken for arbitrary digraphs that may have, in particular, more
than one undominated knot. That is why, the solution obtained for strong digraphs provides no means for
ranking the vertices of an arbitrary digraph.

Having in mind this general case (which was not given much attention in the literature) let us come back to
the system of equations~(\ref{D}). If $\G$ contains more than one undominated knot, there is no spanning
diverging tree in~$\G$. On the other hand, it follows from $\q L=0$ (Theorem~\ref{th4}) that
$L^{\scriptscriptstyle T}\q^{\scriptscriptstyle T}=0,$ that is, every row of $\q$ is a solution to~(\ref{D}).
Then by Proposition~\ref{prop11}, ${\rm rank}\,\q=v$ and ${\rm rank}\,L=n-v$, where $v$ is the forest
dimension of $\G$ (which is the number of undominated knots in~$\G$). Consequently, the system (\ref{D}) has
exactly $v$ linearly independent solutions. If $j\in K_i$ and $|K_i|=k$, then by Theorem~\ref{th2}$'$ the
$j$th row of $\q$ has the form $(0,0\cdc 0,t_{i_1}\cdc t_{i_k},0\cdc 0)$.  Here
$t_{i_s}=\e(\TT^{i_s})/\e(\TT)$ and the vertices ${i_1}\cdc {i_k}$ that belong to $K_i$ are assumed to carry
neighboring numbers. The multiplier $1/\e(\TT)$ is the same for all the solutions. Dividing by it leads us to
the solutions with $t'_{i_s}=\e(\TT^{i_s}),$ i.e., $t'_{i_s}$ is the weight of the set of spanning trees that
diverge from $i_s$ in the undominated knot $K_i$. This provides a description of $v$ linearly independent
solutions to~(\ref{D}).

Thus, by solving (\ref{D}) we do not generally obtain a one-dimensional family of evaluation vectors. For
every undominated knot in $\G$, there is a corresponding partial solution to~(\ref{D}). The general solution
is provided by all their linear combinations. A reasonable ultimate vector of estimates is the convex
combination of the rows of $\q$ with all weights equal to $1/n$ (the arithmetic mean of the rows of $\q$).
This solution is the probability distribution on the set of digraph vertices implemented in the Markov chains
related to $\G,$ provided that the starting distribution on the set of states is uniform.

For the example in Section~\ref{stru}, such a solution to (\ref{D}) is
\[
x\approx (0;0;0.0911;0.1791;0;0.1549;0.1701;0;0;0.0928;0.1701;0;0.1418)^{\scriptscriptstyle T}.
\]
In this solution, as well as in all other solutions, the vertices that are outside of all undominated knots,
are given zero estimates, which is not always reasonable. The estimates based on the matrices $Q(\tau)$
instead of $\q$ do not offer this feature. The problem of their analysis looks meaty.

\section{Forest matrices and the structure of digraphs}
\label{stru}

As was noted in Section~\ref{sec2}, the union $\ktil=\cupo^{v}_{i=1}K_i$ of undominated knots of a digraph is
considered in the theory of decision making as a natural set of alternatives chosen on the base of a binary
relation (digraph) of preferences~\cite{Schwartz}.

Generally speaking, finding the undominated knots and the vertices reachable from each undominated knot is
the first task in discovering the structure of a digraph. According to item~2 of Theorem~\ref{th2}$',$ the
calculation of $\q$ immediately solves this problem. Indeed, the nonzero columns of $\q$ correspond to the
elements of $\ktil$, and the row numbers of the nonzero elements of such a nonzero column index the vertices
reachable from the corresponding undominated knot. In particular, $j\in\ktil$ and $i\in\ktil$ belong to the
same undominated knot if and only if $\q_{ij}\ne0$.

After the appropriate renumbering of the vertices (giving the first numbers to the vertices in $K_1$, the
following numbers to the vertices in $K_2$, and the last numbers to the vertices in $V(\G)\setminus\ktil)$,
we obtain the matrix $\q'$ of the form of
\begin{equation}
\label{blockm}
\q'=\left[\begin{array}{cc}\q'_1 &0\\
                           \q'_2 &0\\
           \end{array}\right],
\end{equation}
where $\q'_1$ is a square block diagonal matrix whose diagonal blocks consist of strictly positive elements
and correspond to the undominated knots, and the nonzero elements $\q'_{ij}$ of $\q'_2$ correspond to pairs
$(i,j)$ such that $i\in V(\G)\setminus\ktil$ and $i$ is reachable from $K(j)$.

To find the positions of the nonzero elements of $\q$ by means of approximate calculations, the following
statement can be used.

\begin{proposition}
\label{approx} Suppose that $\G$ is a digraph all whose arcs have the unit weight. Then the elements of
$(I+\e^2(\FF)L(\G))^{-1}$ that exceed $\e^{-1}(\FF)$ occupy the same positions as the nonzero elements
of~$\vj(\G)$.
\end{proposition}

In essence, Proposition~\ref{approx} is formulated for nonweighted digraphs. This is because the structure of
a digraph does not depend on the weights of its arcs.

Thus, the nonzero elements of $\vj$ can be found as follows.

1. Calculate $Q(\tau)=(I+\tau\,L)^{-1}$ with $\tau=\e^2(\FF)$.

2. Replace with zero all the elements of $Q(\tau)$ less than ${\e^{-1}(\FF)}$. The nonzero elements of $\q$
occupy the complement positions.

If $\q$ is approximately calculated by means of the first item of this algorithm, then the accuracy can be
estimated by the value of the elements that are replaced with zeros (or by $\e^{-1}(\FF)$ in the absence of
such elements).
\smallskip

The reachability matrix of digraph is the matrix $(r\_{ij})$ with elements
\[
r\_{ij}=
\begin{cases}
1, &\text{if $j$ is reachable from $i$,}\\
0, &\text{if $j$ is unreachable from $i$.}
\end{cases}
\]

\begin{proposition}
\label{reach} The reachability matrix can be obtained from the matrix $\q^{\intercal}(\tau)$ with any
$\tau>0$ by the replacement of all its nonzero elements with~$1.$ The result is independent of the weights of
arcs$,$ so they can be set equal to~$1.$
\end{proposition}

This proposition follows from Theorem~$1'$.

An algebraic way to reveal the strong components of a digraph is to find the equal rows (or columns) of the
reachability matrix: their equality means that the corresponding vertices belong to the same strong
component)~\cite{Zykov69}. A variant of this algorithm is to calculate the {\it mutual reachability matrix},
which is the Hadamard (componentwise) product of the reachability matrix and its transpose.

The standard means of finding the reachability matrix of a digraph is calculating $(I+A)^{n-1}$, where $A$ is
the adjacency matrix, or the successive calculation of the power matrices $(I+A)^k$ until the stabilization
of the positions of nonzero elements; in both cases, with the replacement of nonzero elements in the
resulting matrix by ones~\cite{Zykov69}.

\begin{example}
Let us calculate $\q$ for the digraph $\G$ shown in Fig.~2 and use it to reveal the structure of~$\G$.

\addvspace{4mm}
\input rpe_f3e.pic
\vskip-3mm
\nopagebreak
\hspace{68mm} Figure~2

\vskip10mm
\smallskip
{\sloppy The weights of arcs are as follows: $
 \e(2,12)=1.33;\;\e(8,2)=1.5;\;\e(13,8)=0.9;\; \e(11,8)=1.1;
\;\e(7,4)=0.95;\;\e(7,5)=1.3;\;\e(7,9)=1.4; \;\e(5,9)=1.6; \;\e(6,9)=1.25;\;\e(6,3)=1.7;\; \e(3,10)=1.67;\;
  \e(4,13)=1.2;\;\e(13,5)=1.2
$; the weights of the remaining arcs are equal to~one.

}

The result of the approximate calculation of $\q$ by means of Theorem~\ref{Tl} is as follows:

\def\iH{0.5em}\def\iI{0.9em}\def\iJ{3.1em}\def\iK{3.2em}\def\iL{1.0em}
\def\iM{2.7em}\def\iN{3.2em}\def\iO{2.3em}
\def\iP{2.8em}\def\iQ{0.3em}\def\iR{3.8em}\def\iS{1.3em}

\begin{tabbing}
.\hspace{5.8em}\=. \hspace{\iH}\=.\hspace{\iJ}\=.\hspace{\iK}\=.\hspace{\iL}\=.\hspace{\iM}
            \=.\hspace{\iN}\=.\hspace{\iL}\=.\hspace{\iI}\=.\hspace{\iO}
            \=.\hspace{\iK}\=.\hspace{\iI}\=.\hspace{\iJ}\=.
\kill
\>\;1\>\;2\>\,\;3\>\;\,4\>\;\,5\>\;6\>\;7\>\;8\>\;9\>\;10\>\;11\>\;12\>\;13\\
\end{tabbing}
\vspace{-2em}
$$
\arraycolsep=0.35em \vj\approx\left[
\begin{array}{rrrrrrrrrrrrr}
0&0&0.1432&0.1267&0&0.2434&0.1203&0&0&0.1458&0.1203&0&0.1003\\
0&0&     0&0.2709&0&     0&0.2573&0&0&     0&0.2573&0&0.2144\\
0&0&0.2690&     0&0&0.4572&     0&0&0&0.2738&     0&0&     0\\
0&0&     0&0.2709&0&     0&0.2573&0&0&     0&0.2573&0&0.2144\\
0&0&0.0916&0.1786&0&0.1557&0.1697&0&0&0.0932&0.1697&0&0.1414\\
0&0&0.2690&     0&0&0.4572&     0&0&0&0.2738&     0&0&     0\\
0&0&     0&0.2709&0&     0&0.2573&0&0&     0&0.2573&0&0.2144\\
0&0&     0&0.2709&0&     0&0.2573&0&0&     0&0.2573&0&0.2144\\
0&0&0.1432&0.1267&0&0.2434&0.1203&0&0&0.1458&0.1203&0&0.1003\\
0&0&0.2690&     0&0&0.4572&     0&0&0&0.2738&     0&0&     0\\
0&0&     0&0.2709&0&     0&0.2573&0&0&     0&0.2573&0&0.2144\\
0&0&     0&0.2709&0&     0&0.2573&0&0&     0&0.2573&0&0.2144\\
0&0&     0&0.2709&0&     0&0.2573&0&0&     0&0.2573&0&0.2144\\
\end{array}\right]
\begin{array}{r}
1\\2\\3\\4\\5\\6\\7\\8\\9\\10\\11\\12\\13\\
\end{array}
$$

This matrix is easily represented in the form~(\ref{blockm}). The first step can be sorting the nonzero
columns (and thus, classifying the elements of undominated knots) by the positions of nonzero entries; this
reveals two undominated knot: $\{3,6,10\}$ and $\{4,7,11,13\}$. On the second step, the rows corresponding to
the other vertices are classified by the positions of nonzero entries;  we conclude that the vertices in the
strong component $\{1,5,9\}$ are reachable from the both undominated knots, whereas the vertices in the
strong component $\{2,8,12\}$ are reachable from the undominated knot $\{4,7,11,13\}$ only. The resulting
matrix is

\begin{tabbing}
.\hspace{8.1em}\=. \hspace{\iP}\=.\hspace{\iP}\=.\hspace{\iR}\=.\hspace{\iJ}\=.\hspace{\iO}
            \=.\hspace{\iJ}\=.\hspace{\iS}\=.\hspace{\iI}\=.\hspace{\iQ}
            \=.\hspace{\iS}\=.\hspace{\iI}\=.\hspace{\iI}\=.
\kill
\>\;\,3\>\;6\>\;10\>\;4\>\;7\;\>\;11\>\;\,13\>\;\;2\>\;\;8\>\;12\>\;1\>\;\,5\>\;\;9\\
\end{tabbing}
\vspace{-2em}
$$
\arraycolsep=0.35em \vj'\approx\left[
\begin{array}{rrrrrrrrrrrrr}
0.2690 &0.4572 &0.2738 &0      &0      &0      &0      &0&0&0&0&0&0\\
0.2690 &0.4572 &0.2738 &0      &0      &0      &0      &0&0&0&0&0&0\\
0.2690 &0.4572 &0.2738 &0      &0      &0      &0      &0&0&0&0&0&0\\
0      &0      &0      &0.2709 &0.2573 &0.2573 &0.2144 &0&0&0&0&0&0\\
0      &0      &0      &0.2709 &0.2573 &0.2573 &0.2144 &0&0&0&0&0&0\\
0      &0      &0      &0.2709 &0.2573 &0.2573 &0.2144 &0&0&0&0&0&0\\
0      &0      &0      &0.2709 &0.2573 &0.2573 &0.2144 &0&0&0&0&0&0\\
0      &0      &0      &0.2709 &0.2573 &0.2573 &0.2144 &0&0&0&0&0&0\\
0      &0      &0      &0.2709 &0.2573 &0.2573 &0.2144 &0&0&0&0&0&0\\
0      &0      &0      &0.2709 &0.2573 &0.2573 &0.2144 &0&0&0&0&0&0\\
0.1432 &0.2434 &0.1458 &0.1267 &0.1203 &0.1203 &0.1003 &0&0&0&0&0&0\\
0.0916 &0.1557 &0.0932 &0.1786 &0.1697 &0.1697 &0.1414 &0&0&0&0&0&0\\
0.1432 &0.2434 &0.1458 &0.1267 &0.1203 &0.1203 &0.1003 &0&0&0&0&0&0\\
\end{array}
\right]
\begin{array}{r}
3\\6\\10\\4\\7\\11\\13\\2\\8\\12\\1\\5\\9\\
\end{array}
$$

Now the digraph $\G$ can be represented in a more descriptive form (Fig.~3).

\addvspace{12mm}
\input rpe_f4e.pic

\hspace{70mm} Figure~3
\end{example}

\section*{Conclusion}

The set of spanning diverging forests of a digraph and the matrix $\vj$ corresponding to the maximum out
forests have been analyzed. It has been established that the matrix $\vj$ coincides with the matrix of
Ces\`aro limiting probabilities of the Markov chains related to the digraph. Therefore, the matrix
$\q^{\intercal}$ can be considered as the matrix of limiting accessibilities of the digraph. Applications of
the matrices of diverging forests to the problems of revealing the digraph structure and scoring from paired
comparisons are discussed.

\section*{\sl\hfill Appendix}

\prolem{\ref{lem2}}{ is given by contradiction. Suppose that for some $k\in\{\1n-v\}$, no forest $F\in\FF_k$
contains any arc of the form $(z,w)$. Consider an arbitrary $F\in\FF_k$ and any vertex $z$ such that
$(z,w)\in E(\G)$. If $z$ is reachable from $w$ in $F$, then we remove from $F$ any arc of the path from $w$
to $z$ (then $z$ becomes unreachable from $w$) and add the arc $(z,w)$. The resulting digraph belongs to
$\FF_k$. Otherwise, if $z$ is unreachable from $w$ in $F$, then we remove from $F$ any arc and add~$(z,w)$.
Again the resulting digraph belongs to $\FF_k$. }

\propro{\ref{prop1}}{ 1. If a vertex is undominated in $\G$, then it is undominated (and thus it is a root)
in every spanning diverging forest of~$\G$.

2. Assume, on the contrary, that the digraph does not contain circuits, but some dominated vertex $j$ is a
root in a maximum out forest. Then, adding an arbitrary arc directed to $j$, we obtain a forest. This implies
that the previous forest was not maximum. }

\prolem{\ref{lem3}.}{ 1. It follows from item~1 of Proposition~\ref{prop1} that $\FF_k^{t\rightarrow
t}=\FF_k$, therefore, $\e(\FF_k)=\e(\FF_{k}^{t \rightarrow t})$ 2. By Lemma~\ref{lem2}, $\FF_k$ contains at
least one forest $F_k$ that has an arc directed to $w$, therefore, this forest belongs to
$\FF_k\setminus\FF_{k}^{w \rightarrow w}$. Since $\FF_{k}^{t\rightarrow t}=\FF_{k}$ and $\e_{ij}>0$ for all
$(i,j)\in E(\G)$, we have $\e(\FF_{k}^{t \rightarrow t})>\e(\FF_{k}^{w\rightarrow w})$. }

\prolem{\ref{lem4}.}{ Suppose that there is a path from $i$ to $j$ in $\G$, but a maximum out forest $F$ does
not contain such paths. Vertices $i$ and $j$ can belong to (a)~the same tree in $F$ or (b)~different trees.
In the case (a), either $j$ is the root of this tree, and then $i$ is reachable from $j$, or the tree
contains an arc $(k,j)$ with some $k\ne i$, as required. Let us show that in the case (b), $j$ cannot be a
root in~$F$. This will imply that $F$ contains an arc $(k,j)$ with $k\ne i$, and the lemma will be proved.
Assume the contrary and consider an arbitrary path in $\G$ from $i$ to $j$ (such a path does exist by the
hypothesis of the lemma): $(i,i_1), (i_1,i_2)\cdc (i_s,j)$. Let $i_0=i$. By the hypothesis of case~(b), $i_0$
is unreachable from $j$ in~$F$. Let $t\in\{0\cdc s\}$ be the maximum number such that $i_t$ is unreachable
from $j$ in~$F$. By removing in $F$ the arcs directed to $i\_{t+1}\cdc i_s$ (the number of these arcs
is~$s-t$), we obtain a spanning diverging forest where $i\_{t+1}\cdc i_s,j$ are roots. By adding $s-t+1$
arcs, $(i_t,i_{t+1})\cdc (i_s,j)$, to this forest, we obtain a digraph with no semicircuits (by the
definition of $t$) and with all indegrees not exceeding~1. The obtained spanning diverging forest has more
arcs, than the maximum out forest $F$ does.  This contradiction completes the proof.
} 

\prolem{\ref{lem5}.}{ Suppose that $(i,j)\in E(\G)\setminus E(F).$ Vertices $i$ and $j$ can belong to (a) the
same tree in $F$ or (b) different trees. In the case (a), the proof is the same as for Lemma~\ref{lem4}.
Consider the case (b). Since $F$ is maximum, $j$ is not a root in $F$ (otherwise, the addition of $(i,j)$
would produce a forest with a greater number of arcs). Consequently, there exists an arc $(k,j)$ from some
vertex $k\ne i$ to~$j$.

Conversely, if $F$ contains an arc $(k,j)$ with $k\ne i$ or $i$ is reachable from $j$ in $F$, then $(i,j)$ is
not included in~$F$. Indeed, otherwise id$(j)\ge 2$ would hold in the first case, and a circuit would occur
in the second case, but both are impossible in a forest.
} 

\propro{\ref{roots_basis}}{ Let $W$ be the set of roots of some maximum out forest in~$\G$. Then all vertices
in $\G$ are reachable from $W$ by the definition of maximum out forest, and the elements of $W$ are mutually
unreachable by Proposition~\ref{prop2}. Hence, $W$ is a vertex basis of~$\G$.

Conversely, let $W$ be a vertex basis of~$\G$. Let us demonstrate that $W$ is the set of roots of some
maximum out forest in~$\G$. The following statement is obtained in \cite{Fiedler}, and we give its proof here
for the sake of completeness.

\begin{lemma}
\label{lem6} For any strong digraph $\G$ and any vertex $j,$ $\G$ contains a spanning tree diverging
from~$j$.
\end{lemma}

\prolem{\ref{lem6}.}{ We construct the desired tree by the sequential addition of arcs. Let $F_0$ be the
subgraph of $\G$ with $V(F_0)=\{j\}$ and $E(F_0)=\emptyset$. Suppose that a subgraph $F_k$ is already defined
and $(z,i)$ is an arbitrary arc in $\G$ such that $z\in V(F_k)$ and $i\in V(\G)\setminus V(F_k)$. Define
$F_{k+1}$ by setting $V(F_{k+1})=V(F_k)\cup\{i\}$ and $E(F_{k+1})=E(F_k)\cup\{(z,i)\}$. It is obvious that
every subgraph $F_k$ is a tree diverging from $j$ and that the definition process does not terminate until a
spanning tree diverging from $j$ is built. Indeed, if the definition process stops at some $F_k$ with
$k<n-1$, this means that $\G$ contains no paths from $V(F_k)$ to $V(\G)\setminus V(F_k)$, thus $\G$ is not
strong.}

To complete the proof of Proposition~\ref{roots_basis}, observe that, by Proposition~\ref{proZy}, the set $W$
is made up of vertices taken singly from each undominated knot of~$\G$. Using Lemma~\ref{lem6}, one can
construct diverging trees rooted at each a such vertex and spanning in their undominated knots. Uniting the
arc sets of these trees with the set $E(F)\setminus(\ktil\times\ktil)$, where $F$ is an arbitrary maximum out
forest in $\G$, we obtain the arc set of the desired maximum forest. Indeed, the constructed digraph is a
diverging forest with the set $W$ of roots and its number of arcs is no less than in~$F$.
} 

\propro{\ref{prop6a}}{ According to Proposition~\ref{roots_basis}, $K^+_i$ contains only one root of~$F$.
Hence, the restriction of $F$ to $K^+_i$ is a tree diverging from this root, since the vertices in $K^+_i$
are unreachable from the other roots of~$F$.
}

\propro{\ref{prop4}}{ Let us prove item~2; item~1 is proved similarly. Consider an arbitrary forest
$F\in\PP^{K\rightarrow i}$. It follows from the definition of undominated knot that all elements of $K$ are
roots in~$F$. In every tree $T\in\TT^j$, the indegree of every vertex, except the root, is one. Therefore,
$F'=T\cup F$ is a spanning diverging forest of~$\G$. $F'$~has the maximum possible number of arcs, since both
$T$ and $F$ are maximum, i.e., $F'\in\FF^{j\rightarrow i}_{n-v}$.  Consequently, $\TT^j\odot\PP^{K\rightarrow
i}\subseteq \FF^{j\rightarrow i}_{n-v}$. Now consider an arbitrary forest $F'\in\FF^{j\rightarrow i}_{n-v}$.
By Proposition~\ref{prop6a}, the restriction of $F'$ to $K$ is a spanning diverging tree in $K$ rooted
at~$j$. Denote it by $T^j$. Let us demonstrate that the forest whose arcs are the remaining arcs of $F'$
belongs to $\PP^{K\rightarrow i}$. Indeed, $i$ is reachable from $K$ in this forest, and if it is not
maximum, then joining the tree $T^j$ with an arbitrary forest of $\PP$ produces a forest with a greater
number of arcs than in $F'$, contradiction. Thus, $F'\in\TT^j\odot \PP^{K\rightarrow i}$, therefore,
$\FF_{n-v}^{j\rightarrow i}\subseteq \TT^j\odot\PP^{k\rightarrow i}$. Item~2 of Proposition~\ref{prop4} is
proved.
} 

\propro{\ref{algo}}{ Item~1 follows from the fact that the identification of all roots transforms any
diverging forest to a diverging tree, whereas the procedure of root splitting described in step~5a produces a
forest diverging from the vertices that constitute the root of the tree.

Item~2. All subgraphs produced by algorithm~1--6 are maximum out forests, since, by construction, the
indegrees of all vertices, except for $v$ roots lying in $K_1^+\cdc K_v^+$, are equal to~1, and the
constructed subgraphs contain no circuits, as so does~$\G^*$. Finally, this algorithm generates {\it all\/}
the maximum out forests of $\G$, since the restriction of a maximum diverging forest to $K_i^+\,(i=1\cdc v)$
is a diverging tree (Proposition~\ref{prop6a}), and the restriction to $T_i\,(i=1\cdc s)$ is a diverging
forest whose roots are exactly the vertices to which the arcs from outside are directed.
} 

\prothe{\ref{th2}}{ 1. The definition of $Q_{n-v}$ implies $\suml^n_{j=1} q^{n-v}_{ij} =$
$\suml^n_{j=1}\e(\FF^{j\rightarrow i}_{n-v})$. Since $\FF^{j_1\rightarrow i}_{n-v}\capo\FF^{j_2\rightarrow
i}_{n-v}=\emptyset$
whenever $j_1\ne j_2$,
and $\cupo^n_{j=1}\FF^{j\rightarrow i}_{n-v}=\FF_{n-v}$, we obtain
\[
\e(\FF_{n-v})=\e(\cupo^n_{j=1}\FF^{j\rightarrow i}_{n-v})= \suml^n_{j=1}\e(\FF^{j\rightarrow i}_{n-v}).
\]

2. Let $q^{n-v}_{ij}=\e(\FF_{n-v}^{j\rightarrow i})\ne 0$. Then $i$ is reachable from $j$ in $\G$ and, by
Proposition~\ref{roots_basis}, $j\in \ktil$. Let us prove the converse statement. Suppose that $j\in\ktil$
and $i$ is reachable from $j$ in~$\G$. By Proposition~\ref{roots_basis}, $j$ is a root in some maximum out
forest of~$\G$. Denote this forest by~$F$. Suppose that $i$ is unreachable from~$j$ in~$F$. According to
Lemma~\ref{lem4}, $i$ cannot be a root in~$F$. Suppose that the arcs $(j,i_1),\:(i_1,i_2)\:\cdc\:(i_s,i)$
make up a path from $j$ to $i$ in~$\G$. Remove from $E(F)$ all the arcs directed to the vertices $i_1\cdc
i_s$, $i$ (the number of them does not exceed $s+1$) and add the arcs $(j,i_1),\:(i_1,i_2)\:\cdc\:(i_s,i)$.
The resulting subgraph $F'$ is also a maximum out forest, and its tree that contains $i$ is rooted at $j$,
therefore, $q^{n-v}_{ij}\ne 0$. Item~2 is proved.

The first statement of item~3 follows from Proposition~\ref{prop4}, the second statement from
Proposition~\ref{prop6a}.

Item~4 is valid, since, by Proposition~\ref{roots_basis}, the sets $\FF_{n-v}^{j\rightarrow j}\,$ $(j\in K)$
make up a partition of~$\FF_{n-v}$.

5. By item~3, if $j_1,j_2\in K$ and $i$ is reachable from $K$, then
\[
\frac{q^{n-v}_{ij_2}}{q^{n-v}_{ij_1}}= \frac{\e(\TT^{j_2})\e(\PP^{K\rightarrow i})}
{\e(\TT^{j_1})\e(\PP^{K\rightarrow i})}= \frac{\e(\TT^{j_2})}{\e(\TT^{j_1})}.
\]
If $i$ is not reachable from $K$, then $q^{n-v}_{ij_1}= q^{n-v}_{ij_2}=0$. Thereby, the desired equality and
thus Theorem~\ref{th2} are proved.
} 

\prothe{\ref{th3}}{ 1. $\q_{ii}\ge\q_{ji}$, since every vertex $i$ is reachable from itself.

2. If $\q_{ii}>\q_{ji}$, then there exists a maximum out forest $F\in\FF^{i\rightarrow i}_{n-v}
\setminus\FF^{i\rightarrow j}_{n-v}$ where $i$ is a root and $j$ is not reachable from $i$. By
Proposition~\ref{roots_basis}, $i\in\ktil$ and, by item~3 of Theorem~\ref{th2}$',$ $j\notin K^+(i)$. Vertex
$i$ is unreachable from $j$ in $\G$, since otherwise $j\in K(i)$ by the definition of~$K(i)$.

3. In view of item~2, $\q_{ii}>\q_{ji}$ implies $i\in\ktil$ and $j\notin K^+(i)$. Since $\q_{ji}>0$, $j$ is
reachable from $i$, therefore, $j\notin\ktil$. Then, by Proposition~\ref{roots_basis}, $j$ cannot be a root
in any maximum out forest.

4. If $\q_{ij}>0$, then $i$ is reachable from $j$ and, by item~2, $\q_{ii}>\q_{ji}$ is impossible. Then, by
item~1, $\q_{ii}=\q_{ji}$ holds. }

\prothe{\ref{I}}{ Let $\q^{2}=(\q^{(2)}_{ij})$. For any $i,j\in V(\G)$, we have
\begin{equation}
\label{qq} \q^{(2)}_{ij}=\suml^{n}_{k=1}\q_{ik}\q_{kj}.
\end{equation}
Consider the case $\q^{(2)}_{ij}\ne 0$.

If $\q_{jj}>\q_{kj}\ne 0$, then, by item~3 of Theorem~\ref{th3}, $k$ is not a root in any maximum out forest
and $\q_{ik}=0$, hence, $\q_{ik}\q_{kj}=0$. Since for all $k \in V(\G)$, $\q_{jj}\ge\q_{kj}$ holds (by item~1
of Theorem~\ref{th3}), for all nonzero terms in the right-hand side of (\ref{qq}) we have
$\q_{kj}=\q_{jj}>0$, consequently,
\begin{equation}
\label{summa} \q^{(2)}_{ij}=\q_{jj}\suml_{k\in K'}\q_{ik},
\end{equation}
where $K'$ is the set of vertices $k\in V(\G)$ such that $\q_{ik}\q_{kj}\ne0$. Observe that $\q_{ik}\ne0$ and
$\q_{kj}\ne0$ are true together iff $j\in\ktil$, $k\in K(j),$ and $i$ is reachable from $K(j)$ (see item~2 of
Theorem~\ref{th2}$'$). That is why $K'=K(j)$. Using (\ref{summa}) and item~3 of Theorem~\ref{th2}$'$, we
obtain
\[
\q^{(2)}_{ij} =\frac{\e(\TT^j)}{\e(\TT)}\suml_{k\in K(j)} \frac{\e(\TT^k)\e(\TT^{K(j)\rightarrow
i})}{\e(\FF_{n-v})}= \frac{\q_{ij}}{\e(\TT)}\suml_{k\in K(j)}\e(\TT^k)=\q_{ij}.
\]

Suppose now that $\q^{(2)}_{ij}=0$. Then, taking $k=j$ in (\ref{qq}), we conclude that either $j\notin\ktil$
or $j\in\ktil$, but $i$ is unreachable from $K(j)$. By item~2 of Theorem~\ref{th2}$'$, this implies
$\q_{ij}=0$. Theorem~\ref{I} is proved. }

\prothe{\ref{th4}}{ Let us prove the equivalent statement $LQ_{n-v}=Q_{n-v}L=0$. Let $S=(s_{jk})=LQ_{n-v}$.
We will show that $S=0$. By definition, $s_{jk}=\suml^n_{i=1}{\l\_{ji}q^{n-v}_{ik}}=s_1+s_2$, where
$s_1=\suml_{i \neq j}{\l\_{ji}q^{n-v}_{ik}}$ and $s_2=\l\_{jj}q^{n-v}_{jk}$. The number $(-s_1)$ is equal to
the weight of the multiset $\GG^{s_1}$ of weighted 2-digraphs\footnote {A 2-digraph is here a multidigraph
with arc multiplicities not exceeding two. The weight of a 2-digraph is the product of the weights of all its
arcs (including multiple ones). } every element\footnote {The multiset $\GG^{s_1}$ is a set consisting of
pairs $(\h,n_1(\h))$, where $\h$ is a 2-digraph, $n_1(\h)\ge1$ being the multiplicity of $\h$ in~$\GG^{s_1}$.
If $n_1(\h)\ge1,$ then not only $(\h,n_1(\h)),$ but also $\h$ will be called an element of $\GG^{s_1}$; in
this case, we will use the notation $\h\in\GG$. } of which is obtained by the addition of some arc $(i,j)\in
E(\G)$ to some forest from $\FF_{n-v}^{k \rightarrow i}$ ($i=\1n$). The result of this addition is generally
a 2-digraph, because $(i,j)$ can already be in this forest. $\GG^{s_1}$ is a multiset, since this
representation of such a 2-digraph $\h$ is not necessarily unique. In this case, $n_1(\h)$ is the number of
different representations. The weight of $\GG^{s_1}$ is
\[
\e(\GG^{s_1})=\suml\_{\h\in\GG^{s_1}} n_1(\h)\e(\h).
\]

Analogously, $s_2$ is the weight of the multiset $\GG^{s_2}$ of weighted 2-digraphs that consists of pairs
$(\h, n_2(\h))$, $n_2(\h)\ge1,$ whose elements are obtained by the addition of all possible arcs $(i,j)\in
E(\G)$ to all forests from~$\FF_{n-v}^{k\rightarrow j}$. We will prove the equality of $\GG^{s_1}$ and
$\GG^{s_2}$, which will complete the proof of $LQ_{n-v}=0$.

Let us show that $\h\in\GG^{s_1}$ if and only if $\h\in\GG^{s_2}$, and $n_1(\h)=n_2(\h)$.

Suppose that $\h$ is a weighted digraph and $u,w\in V(\h)$. By $\,\h+(u,w)\,$ we denote the 2-digraph with
vertex set $V(\h)$ and the multiset of arcs obtained from $E(\h)$ by the increment of the multiplicity of
$(u,w)$ by~1. Similarly, if $\h$ is a 2-digraph and $u,w\in V(\h)$, denote by $\h'=\h-(u,w)$ the 2-digraph
differing from $\h$ in the multiplicity of arc $(u,w)$ only: $n'((u,w))= \max(n((u,w))-1,0).$

Let $\h\in\GG^{s_1}$. By the definition of $\GG^{s_1}$, $\h=F^{k\rightarrow i}_{n-v}+(i,j)$, where
$F^{k\rightarrow i}_{n-v}\in\FF^{k\rightarrow i}_{n-v}$ for some~$i$. Two cases are possible: (1) $j$ belongs
in $F^{k\rightarrow i}_{n-v}$ to the tree rooted at~$k$ and (2) $j$ does not belong to the tree rooted
at~$k$.

In the case (1), $F^{k\rightarrow i}_{n-v}\in\FF^{k\rightarrow j}_{n-v}$ and thus, $\h=F^{k\rightarrow
i}_{n-v}+(i,j)\in \GG^{s_2}$. In the case (2), $F^{k\rightarrow i}_{n-v}$ does not contain $(i,j)$ and $i$ is
unreachable from~$j$. Consequently, by Lemma~\ref{lem5}, $(t,j)\in E(F^{k\rightarrow i}_{n-v})$ for some
$t\ne i$. Then we obtain $\h-(t,j)=(F^{k\rightarrow i}_{n-v}+ (i,j))-(t,j)\in\FF^{k\rightarrow j}_{n-v}$ and
hence, $\h=(\h-(t,j))+(t,j)\in\GG^{s_2}$.

Suppose now that $\h\in\GG^{s_2}$. Then $\h=F^{k\rightarrow j}_{n-v} +(i,j)$ for some $F^{k\rightarrow
j}_{n-v}\in\FF^{k\rightarrow j}_{n-v}$ and some $i\in V(\G)$, $i\ne j$. Let us show that $\h\in\GG^{s_1}$.
Two cases are possible: (1) $i$ belongs in $F^{k\rightarrow j}_{n-v}$ to the tree rooted at $k$ and (2) $i$
does not belong to the tree rooted at~$k$. In the case (1), $F^{k\rightarrow j}_{n-v}\in\FF^{k\rightarrow
i}_{n-v}$ and, therefore, $\h\in\GG^{s_1}$. In the case (2), $F^{k\rightarrow j}_{n-v}$ does not contain
$(i,j)$ and $i$ is unreachable from~$j$. Consequently, by Lemma~\ref{lem5}, $(t,j)\in E(F^{k\rightarrow
j}_{n-v})$ for some vertex $t\ne i$ such that $t$ belongs in $F^{k\rightarrow j}_{n-v}$ to the tree rooted
at~$k$. Then $\h-(t,j)=(F^{k\rightarrow j}_{n-v} +(i,j))-(t,j)\in\FF^{k\rightarrow t}_{n-v}$ and hence
$\h=\h-(t,j)) +(t,j)\in\GG^{s_1}$.

Let us prove now that for every $\h,$ $n_1(\h)=n_2(\h)$. First, $n_1(\h)$ and $n_2(\h)$ do not exceed~2.
Indeed, by the definitions of $\GG^{s_1}$ and $\GG^{s_2},$ at least three arcs would otherwise be directed to
$j$ and then $\h-(i,j)_1,$ where $(i,j)_1$ is any one of these arcs, would not be a forest. It remains to
prove that $n_1(\h)=2$ iff $n_2(\h)=2$. Indeed, $n_1(\h)=2$ means that there exist $i_1, i_2\in V(\G)$ such
that $i_1\ne i_2$, $E(\h)$ contains $(i_1,j)$ and $(i_2,j),$ $\h-(i_1,j)\in\FF^{k\rightarrow i_1}_{n-v},$ and
$\h-(i_2,j)\in\FF^{k\rightarrow i_2}_{n-v}$. This is {\it equivalent\/} (the proof is below) to the fact that
there exist distinct $i_1\in V(\G)$ and $i_2\in V(\G)$ such that $\{(i_1,j), (i_2,j)\} \subseteq E(\h)$ and
$\{\h-(i_1,j),\h-(i_2,j)\}\subseteq\FF^{k\rightarrow j}_{n-v}$, which, in turn, is equivalent to $n_2(\h)=2$.
To prove the equivalence italicized in the previous sentence, let us formulate the following statement, which
is tantamount to every side of that equivalence:
\begin{equation}
\label{ravno} \h-(i_1,j)\in\FF^{k\rightarrow i_2}_{n-v};\quad \h-(i_2,j)\in\FF^{k\rightarrow i_1}_{n-v}.
\end{equation}
To deduce this from the left member of that equivalence, observe that if, on the contrary,
$\h-(i_1,j)\notin\FF^{k\rightarrow i_2}_{n-v}$, then $(i_1,j)$ belongs in $\h$ to a path from $k$ to $i_2$
and thus, $i_2$ is reachable from $j$ in $\h-(i_1,j)$, which contradicts to the presence of $(i_2,j)$ in
$\h-(i_1,j)$. Similarly, $\h-(i_2,j)\in\FF^{k\rightarrow i_1}_{n-v}$. Further, (\ref{ravno}) immediately
implies the right-hand member of that equivalence. To deduce (\ref{ravno}) from the right-hand member of the
equivalence, observe that $(i_2,j)$ is the unique arc directed to $j$ that belongs to the forest
$\h-(i_1,j)$, and the reachability of $j$ from $k$ in this forest implies the reachability of $i_2$ from $k$.
Consequently, $\h-(i_1,j)\in\FF^{k\rightarrow i_2}_{n-v}$. Similarly, $\h-(i_2,j)\in\FF^{k\rightarrow
i_1}_{n-v}$. Further, $\h-(i_1,j)\in\FF^{k\rightarrow i_2}_{n-v}$ implies that $i_2$ is reachable from $k$ in
$\h$, therefore, $\h-(i_2,j)\in\FF^{k\rightarrow i_2}_{n-v}$. Similarly, $\h-(i_1,j)\in\FF^{k\rightarrow
i_1}_{n-v}$. Thereby, the left-hand member of the equivalence is deduced from~(\ref{ravno}). The identity
$LQ_{n-v}=0$ is proved.

The identity $Q_{n-v}L=0$ is equivalent to the validity of the following equality for all $i,j\in\{\1n\}$:
\begin{equation}
\label{ql} q^{n-v}_{ij}\l\_{jj}=-\suml_{k\neq j} q^{n-v}_{ik}\l\_{kj}.
\end{equation}

The left-hand side of (\ref{ql}) is equal to the weight of the multiset of digraphs obtained by the addition
of all possible arcs $(t,j)$ to all maximum out forests where $i$ belongs to a tree rooted at~$j$. Let this
multiset be $\GG^1$. It is easy to see that $\GG^1$ has no multiple elements. Indeed, every element of
$\GG^1$ is obtained from some forest $F^{j\rightarrow i}_{n-v}$ by the addition of the arc $(t,j)$; both the
$F^{j\rightarrow i}_{n-v}$ and $(t,j)$ are uniquely reconstructed from this digraph.

The right-hand side of (\ref{ql}) is equal to the weight of the multiset $\GG^2$ of digraphs obtained by the
addition of all possible arcs $(j,k)\in E(\G)$ $(k\ne j)$ to all maximum out forests where $i$ belongs to a
tree rooted at~$k$. Let us demonstrate that $\GG^2$ does not contain multiple elements too. Assume, on the
contrary, that some element $\h$ belongs to $\in\GG^2$ with multiplicity greater than one. Then two copies of
$\h$ obtained by the addition of some arcs $(j,k_1)$ and $(j,k_2)$ ($k_1\neq k_2$) to some forests
$F^{k_1\rightarrow i}_{n-v}$ and $F^{k_2\rightarrow i}_{n-v}$, respectively, coincide:
\[
\h =F_{n-v}^{k_1\rightarrow i}+(j, k_1) =F_{n-v}^{k_2\rightarrow i}+(j, k_2).
\]
Under this assumption, $k_1,k_2\in\ktil$ and $(j,k_1), (j,k_2)\in E(\G)$, hence, $j\in\ktil$ and $k_1,k_2\in
K(j)$. Then, by item~3 of Theorem~\ref{th2}, $q_{k_1k_1}^{n-v}=q_{jk_1}^{n-v}$, therefore, $j$ is reachable
from $k_1$ in~$F_{n-v}^{k_1\rightarrow i}$. Consequently, $j$ is reachable from $k_1$ in
$F_{n-v}^{k_1\rightarrow i}-(j,k_2)$ also, hence, $j$ and $k_1$ belong to a circuit in
$F_{n-v}^{k_2\rightarrow i}=(F_{n-v}^{k_1\rightarrow i}-(j,k_2))+(j,k_1)$, which contradicts to the
definition of tree. That is why $\GG^2$ has no multiple elements.

Let us prove the coincidence of $\GG^1$ and $\GG^2$. Let $\h=F^{j\rightarrow i}_{n-v}+(t,j)\in\GG^1$. By
Lemma~\ref{lem4}, $t$ is reachable from $j$ in $F^{j\rightarrow i}_{n-v}$. Consider $z\in\h$ such that
$(j,z)$ is the starting arc of the unique path from $j$ to $t$ in~$\h$. The removal of $(j,z)$ in $\h$
produces a maximum out forest that belongs to~$\FF^{z\rightarrow i}_{n-v}$. Indeed, if the path from $j$ to
$i$ in $F^{j \rightarrow i}_{n-v}$ contains $z$, then the path from $z$ to $i$ is preserved after the removal
of $(j,z)$. Otherwise, if the path from $j$ to $i$ in $F^{j\rightarrow i}_{n-v}$ does not contain $z$, then
the removal of $(j,z)$ preserves this path, and along with the arc $(t,j)$ and the path from $z$ to $t$, it
forms a path from $z$ to $i$ in $\h-(j,z)$. After the addition of $(j,z)$ to the maximum out forest
$\h^1-(j,z)$, we obtain a digraph that belongs to~$\GG^2$.

Let $\h=F^{k\rightarrow i}_{n-v}+(j,k)\in\GG^2$. By removing from $\h$ the last arc, $(t,j)$, of the path
from $k$ to $j$ (which exists by Lemma~\ref{lem5}) we obtain a forest~$F^{j\rightarrow i}_{n-v}$. The
addition of the arc $(t,j)$ to it produces a digraph that belongs to~$\GG^1$. Theorem~\ref{th4} is proved.
} 

\propro{\ref{prop11}}{ Suppose that the vertex set $V(\G)$ is indexed in such a way that the first numbers
are given to the vertices in $K_1$, the following numbers, to the vertices in $K_2$, and so on; the last
numbers are given to the vertices in~$R(\G)$. Then both $L$ and $\vj$ are block lower triangular matrices
with $v+1$ blocks. Every block $L_i$, $i=1\cdc v$, of $L$ coincides with the Kirchhoff matrix of the
restriction, $\G_i$, of $\G$ to~$K_i$. Since $\G_i$ is strong, the matrix-tree theorem for digraphs (see,
e.g.,~\cite{Harary}) and Lemma~\ref{lem6} imply that all minors of order $(|K_i|-1)$ of $L_i$ are strictly
positive. Consequently, the rank of each $i$th diagonal block of $L$ is~$|K_i|-1$. Let us show that the rank
of the last block is equal to its order.

We will use the following notation. Let $\f\subset V(\G)$. Suppose that $L_{-\f}$ is the matrix obtained from
$L$ by deleting the rows and columns corresponding to the vertices in~$\f$; $\G_{(\f)}$ is the multidigraph
obtained from $\G$ by identifying the vertices in $\f$ into the vertex $\f^{\ast}$: every arc in $\G$
incident to a vertex in $\f$ and also to a vertex not in $\f$ induces an arc in $\G_{(\f)}$ incident to
$\f^{\ast}$ and the same second vertex. $\FF_{\f}$ is the set of all spanning diverging forests in $\G$ whose
roots are exactly the vertices in~$\f$. The following statement is due to Fiedler and
Sedl\'a\v{c}ek~\cite{Fiedler}:

\begin{lemma}
\label{lem7} For any $\f\subset V(\G),$ $\det L_{-\f}=\e(\FF_{\f})$.
\end{lemma}

Since the weight of $\FF_{\f}$ is equal to the weight of the set $\TT_{(\f^{\ast})}$ of trees diverging from
$\f^{\ast}$ in $\G_{(\f)}$, we have $L_{-\f}=\e(\TT_{(\f^{\ast})})$.

Let $\f=\ktil$. Then $\TT_{(\f^{\ast})}\ne\emptyset$ and ${\rm det}\,L_{-\f}\ne0$. That is why the rank of
the last block of $L$ is equal to its order and, finally, ${\rm rank}\,L=n-v$.

By virtue of item~2 of Theorem~\ref{th2}$',$ the last (its number is $(v+1)$) block of $\vj$ consists of
zeros. The other blocks are nonzero, and their columns are proportional by item~5 of Theorem~\ref{th2}$'$,
i.e., the rank of every such a block is~1. Therefore, ${\rm rank}\vj\,=v$.
} 

\prothe{\ref{Tl}}{ To prove this fact, it suffices to divide both the numerator
$\suml_{k=0}^{n-v}\,\tau^k\,Q_k$ and the denominator $s(\tau)$ of the formula in Theorem~\ref{th1}$''$ by
$\tau^{n-v}$ and to proceed to the limit as $\tau\to\infty$ using the definition of~$\vj$. }

\propro{\ref{Q1}}{ Suppose that $\alpha\_1=\Big(\suml_{k,t=1}^n\e_{kt}\Big)^{-1}$ and
$P_1=(p_{ij}^1)=$\linebreak $=I-\alpha\_1L(\G)$. Then
\[
p_{ij}^1=
\begin{cases}
\e_{ji}\Big(\suml_{k,t=1}^n\e_{kt}\Big)^{-1}, & j\ne i,
\\
\Big(\suml_{k,t=1}^n\e_{kt}-\suml_{k=1}^n\e_{ki}\Big) \Big(\suml_{k,t=1}^n\e_{kt}\Big)^{-1}, & j=i.
\end{cases}
\]
Thus, $p_{ij}^1$ coincides with the $(i,j)$-entry of the matrix $(\suml_{k,t=1}^n\e_{kt})^{-1} Q_1$ for every
$i,j\in V(\G)$. Now the required statement follows from Proposition~\ref{p7.1} and the obvious inequality
$\alpha\_1<(\maxl_{1\le i\le n}\l\_{ii})^{-1}.$ }

\propro{\ref{p7.2}}{ Let $\left\lceil\,X\,\right\rceil$ be the maximum absolute value of the elements of a
matrix~$X$. Let us take an arbitrary small $\esm>0$ (the designation $\esm$ has nothing in common with the
weights of arcs) and find $k\_0$ such that $\,\left\lceil\,
\frac{1}{k}\,\suml_{t=0}^{k-1}\,P^t-B'\right\rceil<\esm$ for every $k>k\_0$. This will prove the proposition.

Choose $i\_0\in\N$ such that for all $i\ge i\_0$ and $j\in\{0\cdc m-1\},$
\begin{equation}
\label{U1} \left\lceil\,P^{im+j}-B_j\,\right\rceil<\frac{\esm}{2}.
\end{equation}
Set
\begin{equation}
\label{U2} i\_1>\frac{2(i\_0+1)}{\esm}.
\end{equation}
Observe that
\begin{equation}
\label{U3} \suml_{j=0}^{m-1} B'=\suml_{j=0}^{m-1} B_j
\end{equation}
\noindent and that for all $t\in\N$,
\begin{equation}
\label{U4} \left\lceil\,P^t-B'\,\right\rceil\le 1.
\end{equation}
Suppose that $i\_2>i\_1, 0\le j\_2<m$, and $k=i\_2m+j_2$. Then, making use of (\ref{U1})--(\ref{U4}), we
obtain
\begin{eqnarray*}
\left\lceil\,\frac{1}{k}\,\suml_{t=0}^{k-1}\,P^t-B'\,\right\rceil &=&
\left\lceil\,\frac{1}{k}\,\biggl(\suml_{t=0}^{i\_0m-1}(P^t-B')+ \suml_{t=i\_0m}^{i\_2m-1}(P^t-B')+
\suml_{t=i\_2m}^{i\_2m+j\_2} (P^t-B')\biggr)\right\rceil
\\ &\le&
\frac{i\_0m}{k}+\frac{1}{k}
\left\lceil\,\suml_{i=i\_0}^{i\_2-1}\,\suml_{j=0}^{m-1}(P^{im+j}-B')\,\right\rceil
+\frac{j\_2+1}{k}
\\ &\le&
\frac{(i\_0+1)m}{k} +\frac{1}{k}\,\suml_{i=i\_0}^{i\_2-1}\,\suml_{j=0}^{m-1}\,\left\lceil
 P^{im+j}-B_j\,\right\rceil<\frac{i\_0+1}{i\_1}+
\frac{(i\_2-i\_0)m}{k}\,\frac{\esm}{2}<\esm.
\end{eqnarray*}
}
\medskip

{\bf Proof of Corollary~2 from Proposition~\ref{p7.2}.} Item~1. In view of Proposition~\ref{p7.2},
\[
PB= P\frac{1}{m}\suml_{j=0}^{m-1}\lim_{i\to\infty}P^{im+j}=
\frac{1}{m}\suml_{j=0}^{m-1}\lim_{i\to\infty}P^{im+j+1}=
\frac{1}{m}\left(\left(\suml_{j=1}^{m-1}B_j\right)+B_0\right)=B.
\]
Similarly, $BP=B$.

Item~2 of this corollary is just a reformulation of item~1.

Item~3. By (\ref{7.2}) and item~1 of this corollary,
\[
B\,L(\G)=\frac{1}{\alpha}B(I-P)=\frac{1}{\alpha}(B-B)=0.
\]
Similarly, $L(\G)\,B=0.$

Item~4.
\[
B^2= B\left(\frac{1}{m}\suml_{j=0}^{m-1}\lim_{i\to\infty}P^{im+j}\right)=
\left(\frac{1}{m}\suml_{j=0}^{m-1}\lim_{i\to\infty}BP^{im+j}\right)= \frac{1}{m}mB=B.
\]
\medskip 

\propro{\ref{7.p}}{ Item~1 follows from the formula of total probability; item~2 follows from item~1 and
Proposition~\ref{p7.2}. }

\propro{\ref{7.3}}{ Since the spectral radius of $P$ is 1,
\[
\suml_{t=0}^{\infty} (aP)^t=(I-aP)^{-1}
\]
holds true for every $0<a<1$.

Multiplying this identity by $(1-a)$ and making use of the substitution $a=\tau/(\tau+1)$, we obtain
\begin{equation}
\label{Rh} \frac{1}{\tau+1}\,\suml_{t=0}^{\infty}\bigl(\frac{\tau}{\tau+1}P\bigr)^t =(I-\tau(P-I))^{-1}.
\end{equation}
It remains to show that
\[
\lim_{\tau\to\infty}\,\frac{1}{\tau+1} \suml_{t=0}^{\infty}\bigl(\frac{\tau}{\tau+1}P\bigr)^t=B.
\]

Applying the reverse substitution $\tau/(\tau+1)=a$, we shall prove an equivalent (in view of
Proposition~\ref{p7.2}) statement, namely,
\[
\lim_{a\to {1^-}}\,(1-a)\suml_{t=0}^{\infty}(aP)^t=B',
\]
\noindent where $a\to {1^-}$ designates the convergence to 1 from the left. Given a small $\esm>0,$ let us
find $a\_0$ such that for $a\_0<a<1,$
\[
\left\lceil (1-a)\,\suml_{t=0}^{\infty}\,(aP)^t-B'\right\rceil<\esm
\]
holds.

As well as in the proof of Proposition~\ref{p7.2}, take $i\_0\in\N$ such that
\begin{equation}
\label{V1} \left\lceil\,P^{im+j}-B_j\right\rceil<\frac{\esm}{2}
\end{equation}
for all $i\ge i\_0$ and $j\in\{0\cdc m-1\}$.

Choose $a\_0$ such that $0<a\_0<1$ and
\begin{equation}
\label{V2} (1-a\_0)i\_0m<\frac{\esm}{2}.
\end{equation}
Using~(\ref{U3}), (\ref{U4}), (\ref{V1}), (\ref{V2}), and the absolute convergence of the series under
consideration, for every $a\_0 <a<1$ we obtain
\begin{eqnarray*}
\left\lceil\,(1-a)\,\suml_{t=0}^{\infty}\,(aP)^t-B'\,\right\rceil &=&
\left\lceil\,(1-a)\,\suml_{t=0}^{\infty}\,(aP)^t - (1-a)\,\suml_{t=0}^{\infty}\,a^tB'\,\right\rceil
\\ &=&
(1-a)\,\left\lceil\,\suml_{t=0}^{i\_0m-1}\,a^t(P^t-B')+
\suml_{i=i\_0}^{\infty}\,\suml_{j=0}^{m-1}\,a^{im+j}(P^{im+j}-B')\,\right\rceil
\\ &\le&
(1-a)i\_0m+(1-a) \suml_{j=0}^{m-1}\,\suml_{i=i\_0}^{\infty}\,a^{im+j}\,\left\lceil
 P^{im+j}-B_j\,\right\rceil
\\ &<&
\frac{\esm}{2}+ \frac{\esm}{2}(1-a)\,\suml_{j=0}^{m-1}\,\frac{a^{i\_0m+j}}{1-a^m}
=\frac{\esm}{2}+\frac{\esm}{2}\,a^{i\_0m}\le\esm.
\end{eqnarray*}
}

\prothe{\ref{M}}{ The theorem is proved by substituting (\ref{7.1}) in (\ref{Bl}) and applying
Theorem~\ref{Tl} and Proposition~\ref{p7.1}. }

\propro{\ref{bliz}}{ {\it Nonnegativity\/} and the `$\Leftarrow$' part of {\it disconnection condition\/}
follow from Theorem~\ref{th2}$'$; the nonstrict version of {\it diagonal maximality\/} follows from
Theorem~\ref{th3}. By item~3 of Theorem~\ref{th2}$',$  the strict versions of {\it diagonal maximality\/} and
{\it transit property\/} are not fulfilled.

The nonstrict version of {\it transit property\/} is easily proved by contradiction. Assume that for some
$i,k,t\in V(G)$, $\G$ contains a path from $i$ to $k$, $i\ne k\ne t$, and every path from $i$ to $t$ includes
$k$, but $p\_{it}>p\_{ik}$. Then there exists a maximum out forest $F$ where $t$ is reachable from $i$, but
$k$ is unreachable from $i$, which contradicts to the assumption. The first part of item~1 of {\it
monotonicity\/} is satisfied in the nonstrict form by the definition of~$\q$, but is not satisfied in the
strict form by item~2 of Theorem~\ref{th2}$'$. By the same reason, the second part of item~1 of {\it
monotonicity} and the `$\Rightarrow$' part of {\it disconnection condition} are not satisfied.

Let us prove the fulfillment of item~2 of {\it monotonicity} in the nonstrict form. Suppose that the weight
of some arc $(k,t)$ is increased by $\D\e_{kt}$ and the weights of the remaining arcs are preserved.  Denote
the resulting digraph by $\G'$ and set $Q'(\tau)=(I+\tau L(\G'))^{-1}$. Then $\D L=L(\G')-L(\G)=XY$, where
$X=(x_{i1})$ is the column vector with $x_{t1}=1$ and $x_{i1}=0$ for all $i\ne t$, and $Y=(y_{1j}) $ is the
row vector with $y_{1k}=-\D\e_{kt}$, $y_{1t}=\D\e_{kt}$, and $y_{1j}=0$ for all $j\ne k,\;j\ne t$. Since the
matrices $I+\tau L(\G')$ and $I+\tau L(\G)$ are nonsingular and the second one is obtained from the first one
by the addition of $\D L$ with ${\rm rank}\,\D L=1$, by \cite{HoJo} we have
\[
Q'(\tau)=Q(\tau)-\frac{\tau Q(\tau) X Y Q(\tau)}{1+\tau Y Q(\tau) X}=Q(\tau)-\frac{Q(\tau) X Y
Q(\tau)}{\frac{1}{\tau}+Y Q(\tau)X}.
\]
Further, $Q(\tau) X Y Q(\tau)=(a_{ij}(\tau)),$ where $a_{ij}(\tau)=\D\e_{kt}
q_{it}(\tau)(q_{tj}(\tau)-q_{kj}(\tau))$, $\;i,j=\1n$, and $Y Q(\tau)X=\D\e_{kt}(q_{tt}(\tau)-q_{kt}(\tau)).$

We obtain
\begin{equation}
\D\,q_{ij}(\tau)= \frac{\D\e_{kt} q_{it}(\tau)(q_{kj}(\tau)-q_{tj}(\tau))} {\frac{1}{\tau}+\D\e_{kt}
(q_{tt}(\tau)-q_{kt}(\tau))}= \frac{q_{it}(\tau)(q_{kj}(\tau)-q_{tj}(\tau))}
{\frac{1}{\D\,\e_{kt}\tau}+q_{tt}(\tau)-q_{kt}(\tau) },\quad i,j=\1n. \label{deltaqtau}
\end{equation}

Let $Q^{\intercal}(\tau)=P(\tau)=(p_{ij}(\tau))$. Rewrite (\ref{deltaqtau}) for $P(\tau)$:
\begin{equation}
\D\,p_{ji}(\tau)= \frac{p_{ti}(\tau)(p_{jk}(\tau)-p_{jt}(\tau))} { {\frac{1}{\D\,\e_{kt}\tau}}
+p_{tt}(\tau)-p_{tk}(\tau) },\quad i,j=\1n. \label{deltaptau}
\end{equation}

Then for every $i\in V(\G)$,
\[
\D\,p_{it}(\tau)-\D\,p_{ik}(\tau)= \frac{(p_{tt}(\tau)-p_{tk}(\tau))(p_{ik}(\tau)-p_{it}(\tau))} {
{\frac{1}{\D\,\e_{kt}\tau} } +p_{tt}(\tau)-p_{tk}(\tau) }.
\]
Suppose that there exists a path from $i$ to $k$ and every path from $i$ to $t$ contains~$k$. Then, by the
matrix-forest theorem, $p_{ik}(\tau)>p_{it}(\tau)$ and $p_{tt}(\tau)-p_{tk}(\tau)>0$. Proceeding to the limit
as $\tau\to\infty$, we obtain $\D\,p_{it}(\tau)\ge\D\,p_{ik}(\tau)$.

For every vertex not in $\ktil$, the corresponding column of $\vj$ is zero. On the other hand, $\vj$ has no
zero rows, since $\vj$ is stochastic. That is why the {\it reversal property\/} is not satisfied. It is easy
to verify that the {\it triangle inequality for proximities\/} is broken for any $i,j,k\in\{\1n\}$ such that
$i\in\ktil$ and $j,k\in K^+(i)\setminus K(i)$. This implies that {\it metric representability of proximity\/}
is not satisfied either. }

\propro{\ref{approx}}{ By Theorem~\ref{Tl},
\[
\vj=\lim_{\tau\to\infty}\,(I+\tau\,L)^{-1}.
\]

Let us determine $\tau$ such that the calculation of $(I+\tau\,L)^{-1}$ enables one to separate the zero and
nonzero elements of~$\q$. Substituting the notation
\[
A(\tau)=(a_{ij})=\frac{1}{s(\tau)}\suml_{k=0}^{n-v-1}\,\tau^k\,Q_k, \quad
C(\tau)=(c_{ij})=\frac{1}{s(\tau)}\tau^{n-v}\,Q_{n-v}
\]
in (\ref{param}), we obtain
\[
Q(\tau)=A(\tau)+C(\tau).
\]
As $\tau\to\infty$, we have $A(\tau)\rightarrow 0$, $Q(\tau)\rightarrow {Q_{n-v}\slash\e(\FF_{n-v})},$ and
$C(\tau)\rightarrow {Q_{n-v}\slash\e(\FF_{n-v})}.$

Let ${\FF}$ be the set of all spanning diverging forests of~$\G$. By Lemma~2 from~\cite{CheSha97},
$\e(\FF)={\rm det}(I+L).$

Set
\[
\tau=\e^2(\FF)>1.
\]
Then for every $a_{ij}(\tau),$
\begin{equation}
\label{a<} a_{ij}(\tau)<\frac{\tau^{n-v-1} \e(\FF)}{\tau^{n-v}}=\frac{1}{\e(\FF)}
\end{equation}
holds, whereas for any nonzero element $c_{ij}(\tau),$ we have
\begin{equation}
\label{c>} c_{ij}(\tau)>\frac{\tau^{n-v}}{\tau^{n-v}\e(\FF)} =\frac{1}{\e(\FF)}.
\end{equation}

Consequently, all the entries of $Q(\tau)$ that are less than $\e^{-1}(\FF)$, and only such entries,
correspond to the zero entries of~$\vj$.
} 



\begin{thebibliography}{99}

\bibitem{Harary} Harary, F., {\it Graph Theory},
Reading, Mass.: Addison--Wesley, 1969.

\bibitem{Zykov69} Zykov, A.A., {\it Teoriya konechnykh gragov\/}
(Theory of Finite Graphs), Novosibirsk: Nauka, 1969.

\bibitem{Tutte}
Tutte, W.T., {\it Graph Theory}, Reading, Mass.: Addison--Wesley, 1984.

\bibitem{Schwartz} Schwartz, T., {\it The Logic of Collective Choice},
New York: Columbia Univ. Press, 1986.

\bibitem{VIV} Vol'skii, V.I., Choice of Best Alternatives on Directed
Graphs and Tournaments, {\it Avtom. Telemekh.}, 1988, no.~3, pp.~3--17 [{\it Automat. Remote Control},
vol.~49, no.~3, pp.~267--278].

\bibitem{CheSha981} Chebotarev, P.Yu. and Shamis, E.V., On Proximity
Measures for Graph Vertices, {\it Avtom. Telemekh.}, 1998, no.~10, pp.~113--133 [{\it Automat. Remote
Control}, 1998, vol.~59, no.~10, pp.~1443--1459].

\bibitem{CheSha97} Chebotarev, P.Yu. and Shamis, E.V.,
The Matrix-Forest Theorem and Measuring Relations in Small Social Groups, {\it Avtom. Telemekh.}, 1997,
no.~9, pp.~124--136 [{\it Automat. Remote Control}, 1998, vol.~58, no.~9, pp.~1505--1514].

\bibitem{Fiedler} Fiedler, M. and Sedl\'{a}\v{c}ek, J., O
$W$-Bas\'{\i}ch Orientovan\'{y}ch Graf$\ulo$, {\it \v{C}asopis P\v{e}st. Mat.}, 1958, vol.~83, pp.~214--225.

\bibitem{Gelf} Gelfand, I.M., {\it Lektsii po lineynoy algebre\/}
(Lectures on Linear Algebra), Moscow: Nauka, 1971.

\bibitem{Bollo} Bollob\'as, B., {\it Modern Graph Theory}, New York:
Springer, 1998.


\bibitem{SearyRich} Seary, A.J. and Richards, W.D., Partitioning
Networks by Eigenvectors, {\it Proc. Int. Conf. on Social Networks}, London: Univ. of Greenwich Press, 1996,
vol.~1, pp.~47--58.

\bibitem{HoJo} Horn, R.~A. and Johnson C.~R., {\it Matrix Analysis},
Cambridge: Cambridge Univ. Press, 1986.

\bibitem{CheSha982} Chebotarev, P.Yu. and Shamis, E.V., On a Duality
Between Metrics and $\Sigma$-proximities, {\it Avtom. Telemekh.}, 1998, no.~4, pp.~204--209 [{\it Automat.
Remote Control}, 1998, vol.~59, no.~4, pp.~608--612; {\it Erratum}, 1998, vol.~59, no.~10, p.~1501].

\bibitem{David} David, H.A., {\it The Method of Paired Comparisons},
2nd ed., London: Griffin, 1988.

\bibitem{DavidAndr} David, H.A. and Andrews, D.M., Nonparametric
Methods of Ranking from Paired Comparisons, {\it Probability Models and Statistical Analyses for Ranking
Data}, Fligner, M.A. and Verducci, J.S., Eds., New York: Springer, 1993, pp.~20--36.

\bibitem{CookKress} Cook, W.D. and Kress, M., {\it Ordinal Information
and Preference Structures$:$ Decision Models and Applications}, Englewood Cliffs, NJ: Prentice-Hall, 1992.

\bibitem{BelkLev} Belkin, A.R. and Levin, M.S., {\it Prinyatie
resheniy$:$ kombinatornye modeli approksimatsii informatsii\/} (Decision Making: Combinatorial Models of
Information Approximation), Moscow: Nauka, 1990.

\bibitem{CheShaSpr} Chebotarev, P.Yu. and Shamis, E., Constructing an
Objective Function for Aggregating Incomplete Preferences, {\it Econometric Decision Models {\rm (Tangian, A.
and Gruber, J., Eds.).}  Lecture Notes in Economics and Mathematical Systems}, Berlin: Springer,
pp.~100--124.

\bibitem{CheSha99} Chebotarev, P.Yu. and Shamis, E., Preference Fusion
when the Number of Alternatives Exceeds Two: Indirect Scoring Procedures, {\it J. Franklin Inst.}, 1999,
vol.~336, pp.~205--226; Erratum, {\it J. Franklin Inst.}, 1999, vol.~336. pp.~747--748.

\bibitem{Laslier} Laslier, J.-F., {\it Tournament Solutions and
Majority Voting}, Berlin: Springer, 1997.

\bibitem{MoonPullman} Moon, J.W. and Pullman, N.J., On Generalized
Tournament Matrices, {\it SIAM Rev.}, 1970, vol.~12, pp.~384--399.

\bibitem{Berman} Berman, K.A., A Graph Theoretical Approach to
Handicap Ranking of Tournaments and Paired Comparisons, {\it SIAM J. Algebraic Discrete Methods}, 1980,
vol.~1, pp.~359--361.

\bibitem{Maxwell} Maxwell, J.K., {\it Electricity and Magnetism}, 3rd
ed., vol.~1, London: Oxford Univ. Press, 1892.

\end{thebibliography}
\end{document}